\documentclass{amsart}
\usepackage{amssymb, amsmath, amstext, amsopn, amsthm, amscd}
\usepackage{mathtools}
\usepackage{tabularx,ragged2e,booktabs,caption}
\usepackage{tikz}
\usepackage{graphicx}
\usetikzlibrary{positioning}
\usepackage{natbib}
\usepackage{tikz-cd}
\usepackage[linktocpage=true]{hyperref}

\usepackage{algpseudocode}
\usepackage{algorithm}

\newcommand{\Cc}{\mathbb{C}}
\newcommand{\Rr}{\mathbb{R}}
\newcommand{\Nn}{\mathbb{N}}
\newcommand{\Ss}{\mathbb{S}}
\newcommand{\GG}{\mathfrak{g}}
\newcommand{\SL}{\mathfrak{sl}}
\newcommand{\SU}{\mathfrak{su}}
\DeclareMathOperator{\Tr}{Tr}

\newcommand{\trans}{\dagger}
\providecommand{\e}[1]{\ensuremath{\times 10^{#1}}}

\theoremstyle{plain} 
\newtheorem{thm}{Theorem}[]

\theoremstyle{definition}

\newtheorem{rem}{Remark}[] 

\title[Casimir preserving scheme for spherical ideal hydrodynamics]{A Casimir preserving scheme for long-time simulation of spherical ideal hydrodynamics}
\date{}
\author{Klas Modin}
\author{Milo Viviani}
\address{Department of Mathematical Sciences, Chalmers University of Technology and University of Gothenburg, SE-412 96 Gothenburg, Sweden}
\email{klas.modin@chalmers.se, viviani@chalmers.se}

\begin{document}

\begin{abstract}
The incompressible 2D Euler equations on a sphere constitute a fundamental model in hydrodynamics. The long-time behaviour of solutions is largely unknown; statistical mechanics predicts a steady vorticity configuration, but detailed numerical results in the literature contradict this theory, yielding instead persistent unsteadiness. Such numerical results were obtained using artificial hyperviscosity to account for the cascade of enstrophy into smaller scales.  Hyperviscosity, however, destroys the underlying geometry of the phase flow (such as conservation of Casimir functions), and therefore might affect the qualitative long-time behaviour. Here we develop an efficient numerical method for long-time simulations that preserve the geometric features of the exact flow, in particular conservation of Casimirs. Long-time simulations on a non-rotating sphere then reveal three possible outcomes for generic initial conditions: the formation of either 2, 3, or 4 coherent vortex structures. These numerical results contradict the statistical mechanics theory and show that previous numerical results, suggesting 4 coherent vortex structures as the generic behaviour, display only a special case. Through integrability theory for point vortex dynamics on the sphere we present a theoretical model which describes the mechanism by which the three observed regimes appear. We show that there is a correlation between a first integral $\gamma$ (the ratio of total angular momentum and the square root of enstrophy) and the long-time behaviour: $\gamma$ small, intermediate, and large yields most likely 4, 3, or 2 coherent vortex formations. Our findings thus suggest that the likely long-time behaviour can be predicted from the first integral $\gamma$.

\end{abstract}

\maketitle

\tableofcontents
 
\section{Introduction}

The motion of an ideal fluid restricted to the surface of a sphere is a fundamental model in oceanography, meteorology, and astrophysics \citep[see][and references therein]{MaBe2002,Do2012,Pe2013,Ze2018}.
The equations of motion, first studied by Euler in 1757, encode a rich geometry---a \emph{Lie--Poisson structure}---which results in conservation of energy, momentum, and Casimir functions~\citep[see][]{arn,maz,ArKh1998}.

The ultimate `fate' of 2D fluid motion in a bounded domain is largely unknown~\citep{Ne2016}.
Statistical mechanics theories, such as developed by \citet{Mi1990} and \citet{RoSo1991}, are based on maximizing entropy of a coarse-grain probability distribution of the macroscopic states under conservation of energy and (at least some of) the Casimirs.
Such models predict a steady equilibrium of large-scale coherent vortex structures, with a functional relation between vorticity and stream function.

To test the statistical model of Miller, Robert, \&\ Sommeria (MRS) a natural approach is to use long-time numerical simulations.
A serious complication is the `inverse energy cascade' where energy from small scales are eventually fed into large scales whereas enstrophy cascades in the forward direction towards smaller scales.
This process was first described by \citet{Kr1967}. 
Of course, in a numerical simulation the spatial resolution is finite, so one can never fully resolve the fine-scale structure.
As a remedy, a common approach is to adopt a \emph{subgrid model}, most often \emph{hyperviscosity}, to account for the enstrophy cascade to smaller scales \citep[see][and references therein]{QiMa2014}.
The inverse energy cascade is related to the conservation of Casimirs, although the exact relation is unknown.
In addition to energy, circulation (linear Casimir), and enstrophy (quadratic Casimir), there are several numerical investigations reporting that cubic and possibly higher order Casimirs also play a role in the formation of large scale coherent vortex structures~\citep{AbMa2003,DuFr2010}.
On the non-rotating sphere, \citet*{dri} provided numerical evidence that, for randomly generated initial data, 
the long-time behaviour results in a non-steady interaction largely between two positive and two negative coherent vortex structures (referred to as \emph{vortex blobs} in this paper) essentially governed by finite dimensional point vortex dynamics.
Seemingly persistent unsteadiness in numerical solutions of 2D Euler fluids were also reported by \citet{SeKi1998} but for special initial conditions.
\citeauthor{dri} (DQM) argue that, in fact, the unsteady four vortex blobs behaviour is generic.
This statement is in stark contrast to the previous notion that a steady equilibrium is the generic behaviour.
However, DQM used methods with hyperviscosity and in their simulations the percentage decay in enstrophy is between 30-60\%, so hyperviscosity clearly comes into play, but precisely how and if it affects the long-time result is unclear.

In this paper, based on a new numerical method that exactly conserves all Casimir functions thereby eliminating the need for hyperviscosity, we give strong evidence that neither MRS nor DQM are correct. 
Or, in a way, they are both correct---it all depends on the regime of the initial conditions.
Based on the non-dimensional non-negative number $\gamma$ given by the quotient between the total angular momentum and the total enstrophy, we identify three different regimes: generically\footnote{The exact definition of `generic' here is that the initial vorticity is sampled as a random field in the space of $L^2$ functions, as described in \autoref{sec:NonZeroMomsimulation} below.}, if $\gamma \lesssim 0.15$ then most likely 4 vortex blobs form (the behaviour observed by DQM), if $\gamma\gtrsim 0.4$ then most likely 2 vortex blobs form (the behaviour suggested by MRS), and if $0.15\lesssim\gamma\lesssim 0.4$ we have found a new, intermediate regime where most likely 3 vortex blobs form. 
The 2 vortex blobs formation is steady (or at least almost steady), whereas the 3 and 4 blobs formations are unsteady.
Furthermore, through point vortex dynamics, we suggest a theoretical mechanism which qualitatively explains the three regimes.
This theory, which also predicts results observed on the torus, is not based on statistical mechanics (i.e., maximizing entropy, like MRS) but rather on integrability theory (results on quasi-periodicity) for point vortex dynamics.


As mentioned already, the central tool in the discovery of the three regimes is a new numerical scheme for ideal fluids on rotating or non-rotating spheres that encapsulate the full Lie--Poisson geometry (in particular conservation of associated Casimirs).\footnote{It is clear from the definition of $\gamma$ that a scheme with hyperviscosity, such as those used by \citet{dri} with 30-60\% decay in enstrophy but no decay in total angular momentum, could never be used to correctly identify the three regimes.}
It is based on geometric quantization theory developed by \citet{hopPhD,hopyau} in conjunction with the Lie--Poisson preserving numerical time discretization developed by~\citet{modviv1}.
The method can be seen as a spherical analogue of the spatial discretization of the Euler equations on the torus suggested by \citet{ze1} and the associated numerical time discretization suggested by \citet{mcl}.


We now continue the introduction with a more detailed exposition of the equations of motion, an overview of the space and time discretization, and a summary of our main findings.

Consider a homogeneous, incompressible, inviscid, two-dimensional fluid, constrained to the unit sphere $\Ss^2$ embedded in Euclidean $\mathbb{R}^3$, and possibly rotating with constant angular speed about a fixed axis. 
The equations of motion are given by Euler's equations of hydrodynamics
\begin{equation}\label{euleq1}
\begin{array}{ll} 
	\dot{v} + v\cdot\nabla v = -\nabla p -2\widetilde{\Omega} \times v \\
	\nabla\cdot v = 0,
			\end{array}
\end{equation} 
where $v$ is the velocity vector field of the fluid, $p$ is its internal pressure, and $\widetilde{\Omega}=(\Omega\cdot n)n$ is the projection of the angular rotation vector $\Omega\in\mathbb{R}^3$ to the normal~$n$. 
The term $-2\widetilde{\Omega} \times v$ is due to the Coriolis force. 
Equivalent to \eqref{euleq1} is the \textit{barotropic vorticity equation} (also called the \emph{quasi-geostrophic equation} in the case $\widetilde\Omega\neq 0$), formulated in terms of the vorticity variable $\omega=(\nabla\times v)\cdot n$. 
By Stokes' theorem we necessarily have $\int \omega =0$ corresponding to zero circulation. 
Euler's equations~\eqref{euleq1} can now be written
 \begin{equation}\label{euleq2}
\begin{array}{ll}
&\dot{\omega}=\lbrace\psi,\omega\rbrace \vspace{2mm}\\ 
&\Delta\psi = \omega - f,
\end{array}
\end{equation} 
where $f\coloneqq 2\Omega\cdot n$ is the \emph{Coriolis parameter}, $\Delta$ is the Laplace--Beltrami operator, $\{\cdot,\cdot\}$ is the Poisson bracket, 
and the stream function $\psi$ is unique by the additional condition 
\begin{equation*}\label{eq:zero_condition}
	\int \psi =0 .	
\end{equation*}
The vorticity equation \eqref{euleq2} constitutes an infinite dimensional \emph{Lie--Poisson system} \citep[cf.][]{ArKh1998} on the space of smooth zero mean functions
\begin{equation*}	C^\infty_0(\Ss^2) =\Big\lbrace\omega \in C^\infty(\Ss^2)\mid \int\omega = 0 \Big\rbrace.
\end{equation*}
The Hamiltonian is
 \begin{equation}\label{eq:ham_euler}
 H(\omega) = \frac{1}{2}\int (\omega-f)\psi,
 \end{equation}
and the (infinitely many) Casimir functions are given, for any smooth real function $g\in C^\infty(\Rr)$, by 
\begin{equation}\label{eq:casimirs}
	\mathcal C(\omega)=\int g(\omega).
\end{equation}
Often $g$ is chosen as monomials, and the corresponding Casimirs
\begin{equation*}
	\mathcal C_k(\omega)=\int \omega^k, \qquad k=1,2,\ldots
\end{equation*}
are called linear, quadratic, cubic, etc.
Each Casimir \eqref{eq:casimirs} is indeed a first integral:
 \[
  \frac{d}{dt}\int g(\omega)=-\int g'(\omega)v\cdot\nabla \omega=-\int v\cdot\nabla g(\omega)=\int (\nabla\cdot v) g(\omega)=0,
  \]
where we have used that  
  \[\lbrace\psi,\omega\rbrace_p = p\cdot(\nabla\psi_p\times\nabla\omega_p)=(p\times\nabla\psi_p)\cdot\nabla\omega_p=-v_p\cdot\nabla\omega_p
  \]
for any $\psi,\omega\in C^\infty(\Ss^2)$ and any $p\in\Ss^2$.
Notice, in particular, that the Casimirs are conserved for any choice of Hamiltonian; this reflect the underlying Lie--Poisson geometry which is foliated in \emph{co-adjoint orbits} preserved by any Hamiltonian flow \citep[cf.][Ch.~13-14]{MaRa1999}.

The traditional approach to numerical discretization of PDE is to construct schemes of high local order of accuracy, using for example finite element or finite volume schemes. 
Rather than focusing on local accuracy, we take here conservation of the Casimir functions~\eqref{eq:casimirs} and the underlying geometric structure as a guiding principle for spatial discretization: we wish to replace the infinite dimensional Lie--Poisson structure $(C^\infty_0(\Ss^2),\{\cdot,\cdot\})$ by a finite dimensional analogue.
We require the number of conserved quantities to increase with the size of the spatial discretization. 
This cannot be achieved by a truncated spectral decomposition of the vorticity, essentially because the space spanned by a truncated spectral basis is not closed under the Poisson bracket.
Instead, we take the approach proposed by \citet{ze2} based on the theory of geometric quantization studied in \citep{bhss,bms,hopPhD}. 
It provides a sequence $N=1,2,\dots$ of finite dimensional Lie algebras, that converges
to the infinite dimensional Lie algebra of smooth functions on the sphere as $N\to\infty$. 
The sequence is given explicitly by the Lie algebra $\SU(N)$ (or $\SL(N,\Cc)$)\footnote{$\SU(N)$ is the Lie algebra of $N\times N$ skew-hermitian complex matrices with trace zero, $\SL(N,\Cc)$ is the Lie algebra of $N\times N$ complex matrices with trace zero.} for $N=1,2,\ldots$. 
For any choice of $N$ we get an ODE which is a finite dimensional analogue of \eqref{euleq2}
 \begin{equation}\label{euleq3}
\begin{array}{ll}
&\dot{W}=[\Delta^{-1}_N(W-F),W]_N,
\end{array}
\end{equation} 
where $W\in\SU(N)$ (corresponding to the vorticity $\omega$), $F\in\SU(N)$ (corresponding to the Coriolis parameter $f$), $\Delta_N\colon \SU(N)\to \SU(N)$ is the discrete Laplace-Beltrami operator (corresponding to $\Delta$), and $[\cdot,\cdot]_N$ is the rescaled matrix commutator (corresponding to $\{\cdot,\cdot\}$).
The matrix differential equation \eqref{euleq3} is an \emph{isospectral flow}, meaning that the eigenvalues of $W$ are invariant in time.
The conservation of these eigenvalues corresponds to the conservation of the Casimirs.
Exactly how $W$ in \eqref{euleq3} approximates $\omega$ in \eqref{euleq2} is described in a complicated (but explicit) linear change of coordinates between $W$ and a truncated spherical harmonics basis.
Details are given in \autoref{sec:spatial_discretization}.
A feature of the spatial discretization is that $W\mapsto \Delta^{-1}_N (W-F)$ can be computed in only $\mathcal O(N^2)$ operations.
Thus, the main computational complexity is due to matrix multiplications in the bracket $[\cdot,\cdot]$ (which has complexity $\mathcal O(N^3)$).
Details on the computational complexity are given in \autoref{sub:complexity}.
 
To discretize \eqref{euleq3} in time we apply a Lie--Poisson preserving \emph{isospectral symplectic Runge--Kutta} integrator \citep{modviv1}. 
These numerical methods exactly conserve (i.e.\ up to rounding errors) the discrete Casimirs (eigenvalues), they nearly conserve the Hamiltonian (`nearly' in the sense of \emph{backward error analysis of symplectic integrators}, c.f.~\citet{HaLuWa2006}), and they exactly conserve the Lie--Poisson flow structure (in short, this means that the time discretized system correspond to a continuous Lie--Poisson flow on $\SU(N)$ for a slightly modified Hamiltonian).
The IsoSRK integrators are necessarily implicit, thus requiring nonlinear root-solving at each time step.
As a comparison, we also employ the standard explicit Heun method for time discretization of \eqref{euleq3}.



In \autoref{sec:simulations} we present numerical simulations on a non-rotating sphere ($F=0$).
First, in \autoref{sec:DQMsimulation}, we use the same randomly generated initial data as suggested by~\citet{dri}.
Long-time simulations are carried out for both types of time discretizations (IsoSRK and Heun) and various levels of spatial discretization.
Our numerical results verify, but now without hyperviscosity, the formation of a quadruple of vortex blobs moving quasi-periodically with no sign of reaching steadiness.
However, although the DQM initial conditions were randomly generated, we claim they cannot represent the generic behaviour because the total angular momentum is zero.
The motivation by \citeauthor{dri} to set momentum to zero was ``to avoid starting with a flow organised at the largest possible scale''.
Herein lies the implicit assumption that the value of the momentum does not affect the qualitative behaviour.
On the doubly periodic square (i.e., the flat torus) the assumption is correct: momentum does not influence the dynamics and can therefore safely be set to zero.
On the sphere, however, the momentum \emph{strongly} affects the dynamics.
In fact, our results suggest that the generic qualitative behaviour on a non-rotating sphere is essentially \emph{governed} by the value of the total angular momentum.
Indeed, in \autoref{sec:NonZeroMomsimulation} we generate 16 sets of initial vorticity as samples from a Gaussian random field on the space of $L^2$-functions.
In the corresponding 16 long-time simulations we observe the following qualitative behaviour: 5 of them give 4 vortex blobs, 9 of them give 3 vortex blobs, and 2 of them give 2 vortex blobs.
We also observe that the non-dimensional number $\gamma = \| \mathbf{L}\|/(R\sqrt{\mathcal{C}_2})$ (total angular momentum divided by the radius of the sphere times the square root of enstrophy) gives a probabilistic indication of which `qualitative regime' the fluid configuration develops into: small values (approximately less than $0.15$) result in 4 vortex blobs, large values (approximately larger than $0.4$) result in 2 vortex blobs, and intermediate values result in 3 vortex blobs.
The number $\gamma$, computable from the initial conditions, is thus implicated in predicting the fluid's long-time qualitative behaviour.
Of the three regimes, only the 2 vortex formation is steady (up to a constant speed rotation about the momentum axis).

It is natural to ask for a theoretical model explaining the three observed regimes.
Clearly, the statistical mechanics based MRS theory is insufficient; it incorrectly predicts steadiness and does not predict or offer insights to why there should be three regimes.
Instead, we have found a different theory which explains the mechanisms by which the regimes appear: it is closely related to integrability theory for point vortex dynamics (PVD).
Recall that a Hamiltonian system is called \emph{integrable} if there is a local change of variables in which the dynamics is described by quasi-periodic linear motion on tori.\footnote{Equivalently, integrability of a $2n$-dimensional Hamiltonian system can be characterized by the existence of $n$ first integrals in involution \citep[cf.][]{Ar1989}.}\
PVD constitute a class of Hamiltonian $N$-particle systems that describe, at least formally, special solutions to the Euler equations \eqref{euleq1} in the non-rotating case ($\Omega = 0$).
\citet{are} refers to PVD as ``a classical mathematics playground'': although the connection to fluid mechanics has always remained in the background, mathematicians have studied these finite dimensional Hamiltonian systems in their own right, observing that ``many strands of classical mathematical physics come together''~\citep[Sect.~I]{are}.
A frequently addressed question is whether a particular number of point vortices on some given geometry (for example the sphere) yields integrable dynamics or not.
In \autoref{sec:pvd} of this paper we (re)connect the mathematical theory for integrability of PVD to the long-time behaviour of a continuous, generic incompressible fluid, thereby obtaining an explanation of the three observed regimes.
This is briefly how the mechanism works:
\begin{enumerate}
	\item Smaller vortex formations of the same sign merge to larger formations when their trajectories come close enough (the inverse energy cascade).
	\item The motion of $N$ vortex blobs is accurately described by $N$ point vortices as long as the blobs are well-separated (so that no merging occurs).
	A careful, numerical evaluation of this assumption is given in \autoref{sec:pvd_zero_mom}.
	\item If the motion of $N$ vortex blobs is not integrable, then, sooner or later, two vortex blobs of equal sign will reach a point in phase space where they are close enough to merge. 
	\item If, however, the motion of the $N$ vortex blobs is integrable\footnote{Or at least close enough to integrable in the KAM sense, see \autoref{sec:pvd_non_zero_mom}.}\ then the motion remains quasi-periodic with well-separated trajectories and no further merging occurs (integrability acts as a barrier in phase space, preventing further merging of blobs).\footnote{From a mathematical viewpoint, the integrability prevents the dynamical system from being \emph{ergodic}. Ergodicity is assumed in statistical mechanics theories such as MRS.}
\end{enumerate}
To summarize, vortex blobs of equal sign continue to merge until integrability blocks them from doing so.
Thus, in order to find the long-time behaviour, one has to find the largest possible number of point vortices for which the dynamics is integrable.
Here is the key-point: on the non-rotating sphere integrability depends on the total angular momentum.
A 4-PVD system on the sphere is integrable if the momentum is zero, but non-integrable if the momentum is non-zero~\citep{Sa2007}.
If momentum is close to zero one still obtains `integrable like' dynamics since integrable systems are stable in the sense of KAM theory for small perturbations (the small momentum configuration can be viewed as a perturbation of a zero momentum configuration). 
This explains why 4 vortex blobs is the stable long-time regime for fluid configurations with a small $\gamma$ parameter.
If the momentum in a 4 blobs configuration is above the threshold where KAM can be applied, the dynamics is chaotic and sooner or later two of the blobs will merge into a 3 blobs configuration.
Since 3-PVD systems on the sphere are integrable (regardless of the momentum), this explains the intermediate 3 blobs regime.
It remains to explain why 2 blobs are sometimes formed.
If $\gamma$ is large enough, there are already two dominant vortex blobs from the start, so the smaller vortex formations are directly merged with these two without passing through the stable 3 vortex blobs regime.
We thereby have an explanation of the mechanism leading to the three observed regimes.

Conclusions and an outlook to future research are presented in \autoref{sec:conclusions}.
Although our main focus is with the non-rotating sphere, we have included in Appendix~\ref{sec:rossby} numerical examples of Rossby--Haurwitz waves on a rotating sphere, to illustrate the usability of the new method also in the rotating case (relevant for quasi-geostrophic flows in atmospheric dynamics).

\medskip

\noindent\textbf{Acknowledgements.} 
The work was supported by EU Horizon 2020 grant No 691070, by the Swedish Foundation for International Cooperation in Research and Higher Eduction (STINT) grant No PT2014-5823, by the Swedish Foundation for Strategic Research grant ICA12-0052, and by the
Swedish Research Council (VR) grant No 2017-05040.
The authors would also like to thank D.~Dritschel for providing us with the code `Hydra'.

\section{Numerical integration algorithm}\label{sec:num_scheme}

For spatial discretization we use the system of differential equations developed by \citet{ze2}, based on the work of Hoppe et al.\ on the approximation of infinite dimensional Lie algebras \citep{bms,bhss}. 
The Poisson algebra of smooth functions on the sphere is approximated by the finite dimensional matrix Lie algebras $\SL(N,\Cc)$, for the Poisson algebra $C_0^\infty(\Ss^2,\Cc)$, and $\SU(N)$ for the Poisson algebra $C_0^\infty(\Ss^2,\Rr)$.
To discretize the equations in time we use the class of isospectral symplectic Runge--Kutta methods developed by \citet{modviv1}.

\subsection{Spatial discretization via geometric quantization}\label{sec:spatial_discretization}

This section is devoted to the technique used to get a finite dimension analogue of the Euler equations on a sphere. 
The main theoretical concept behind the approach is the so called $L_\alpha$-approximation.   

\subsubsection{$L_\alpha$-approximation}
Consider a Lie algebra $(\GG,[\cdot,\cdot])$ and a family of labeled Lie algebras $(\GG_\alpha,[\cdot,\cdot]_\alpha)_{\alpha\in I}$, where $\alpha\in I=\Nn$ or $\Rr$. 
Furthermore, assume that to any element of this family, a distance $d_\alpha$ and a surjective projection map $p_\alpha:\GG\rightarrow\GG_\alpha$ are associated. 
Then an $L_\alpha$-approximation $(\GG_\alpha,[\cdot,\cdot]_\alpha)_{\alpha\in I}$ of $(\GG,[\cdot,\cdot])$ should fulfill:
\begin{enumerate}
\item if $x,y\in\GG$ and $d_\alpha(p_\alpha(x),p_\alpha(y))\rightarrow 0$ as $\alpha\rightarrow\infty$, then $x=y$;
\item for all $x,y\in\GG$ we have $d_\alpha(p_\alpha([x,y]),[p_\alpha(x),p_\alpha(y)]_\alpha)\rightarrow 0$ as $\alpha\rightarrow\infty$;
\item for $\alpha\gg0$ the projections $p_\alpha$ are surjective.
\end{enumerate}
The above definition is given in \citep{bms}; it is a weak requirement to obtain a limit for a sequence of Lie algebras. 

Consider now the smooth complex functions on the sphere with vanishing mean, denoted $C^\infty_0(\Ss^2,\Cc)$. 
This vector space is endowed with a Poisson structure $\lbrace\cdot,\cdot\rbrace$ given by the skew symmetric bilinear form on $C^\infty_0(\Ss^2,\Cc)$ 
\begin{equation}\label{bracket1}
\lbrace f,g\rbrace(x)=\lvert X_f(x),X_g(x),x \rvert,
\end{equation}
where $X_h(x)=x\times\nabla h(x)$ is the Hamiltonian vector field associated with the Hamiltonian function $h\in C^\infty_0(\Ss^2,\Cc)$.
With this bracket, $C^\infty_0(\Ss^2,\Cc)$ becomes an infinite dimensional Poisson algebra; in particular, it is an infinite dimensional Lie algebra.

A basis for $C^\infty_0(\Ss^2,\Cc)$ is given by the complex spherical harmonics, expressed in the standard azimuthal-inclination coordinates $(\phi,\theta)$ by
\[
Y_{lm}(\phi,\theta)=\sqrt{\dfrac{2l+1}{4\pi}\dfrac{(l-m)!}{(l+m)!}}P^m_l(\cos\theta)e^{\mathrm{i}m\phi},
\quad l\geq 1, \quad m=-l,\dots,l ,
\]
where $P^m_l$ are the \emph{associated Legendre polynomials} (i.e.\ solutions to the general Legendre equation).
Using this basis, an explicit approximating sequence for $C^\infty_0(\Ss^2,\Cc)$ was constructed by \citet{hopPhD}.
The sequence is given by the matrix Lie algebras $(\SL(N,\Cc),[\cdot,\cdot]_N)_{N\in\Nn}$, where $[\cdot,\cdot]_N \coloneqq N^{3/2}[\cdot,\cdot]$ is a rescaling of the matrix commutator $[\cdot,\cdot]$.
The distances $d_N$ are given by suitable matrix norms, and the projections $p_N$ are obtained by associating to each spherical harmonic $Y_{lm}$ a matrix $\mathrm{i}T^N_{lm}\in\SL(N,\Cc)$ defined by
\[
(T^N_{lm})_{m_1m_2}=(-1)^{[\frac{N-1}{2}]-m_1}\sqrt{2l+1} 
\left( \begin{array}{ccc}
\frac{N-1}{2} & l & \frac{N-1}{2} \\
-m_1 & m & m_2 
\end{array} \right),
\]
where the bracket denotes the Wigner 3j-symbols. 
The following $L_\alpha$-convergence result for this approximating sequence have been established:

\begin{thm}[\citet{bhss,bms}]\label{thm:Lalpha}
Consider the Poisson algebra $(C^\infty_0(\Ss^2,\Cc),\lbrace\cdot,\cdot\rbrace)$ with Poisson bracket defined by \eqref{bracket1}. 
Then, for the projections $p_N$ and any choice of matrix norms $d_N$, $(\SL(N,\Cc),[\cdot,\cdot]_N)_{N\in\Nn}$ is an $L_\alpha$-approximation of $(C^\infty_0(\Ss^2,\Cc),\lbrace\cdot,\cdot\rbrace)$.
\end{thm}

\subsubsection{The quantized system}

We can now derive the spatial discretization of the Euler equations via the $L_\alpha$-approximation in Theorem~\eqref{thm:Lalpha}, thereby obtaining a finite dimensional `quantized' system.
We begin without the Coriolis parameter.

For any $N\in\Nn$ an analogue of the Euler equations \eqref{euleq2} is the following flow of matrices
\begin{equation}\label{quant_euler}
\dot{W}=[\Delta^{-1}_N W,W]_N,
\end{equation}
where $W\in\SL(N,\Cc)$ and $\Delta^{-1}_N$ is the inverse of the \emph{discrete Laplacian}, given by the following formula of \citet{hopyau}
\begin{equation}\label{eq:discrete_laplacian}
\Delta_N = \frac{N^2-1}{2}\left([X^N_3,[X^N_3,\cdot]] - \frac{1}{2}[X^N_{+},[X^N_{-},\cdot]] - \frac{1}{2}[X^N_{-},[X^N_{+},\cdot]]\right),
\end{equation}
where $X^N_{\pm}\varpropto T^N_{1\pm 1}$, $X^N_3\varpropto T^N_{10}$. 
The crucial property of $\Delta^{-1}_N$ is that $\Delta^{-1}_NT^N_{lm}=(-l(l+1))^{-1}T^N_{lm}$, for any $l=1,\ldots,N$, $m=-l,\ldots,l$.
That is, the basis elements $T^N_{lm}$ are eigenvectors of the discrete Laplacian $\Delta_N$, which is a direct analogue to the continuous case where the spherical harmonics $Y_{lm}$ are eigenvectors of the Laplace--Beltrami operator $\Delta$.

Let us again, now explicitly, discuss the connection between the continuous vorticity equation \eqref{euleq2} and the quantized version \eqref{quant_euler}.
First, notice that \eqref{quant_euler} is an \emph{isospectral flow}; it preserves the eigenvalues of $W=W(t)$.
This isospectral property is a direct analogue of preservation of Casimirs in the continuous flow \eqref{euleq2}.
Given a continuous vorticity function expanded in the spherical harmonics basis, $\omega = \sum \omega^{lm} Y_{lm}$, the projection operator $p_N$ is given by
\begin{equation*}
	p_N(\omega) = \sum_{l=1}^{N-1}\sum_{m=-l}^l \mathrm{i}\omega^{lm}T^N_{lm} .
\end{equation*}
If the continuous vorticity $\omega$ is real valued, then the spherical harmonics coefficients fulfill $\omega^{lm} = (-1)^m \omega^{l(-m)}$.
The corresponding condition on the matrix $W\in \SL(N)$ is $W + W^\dagger = 0$, i.e.\ it belongs to the subalgebra $\SU(N)$ of trace-free skew Hermitian matrices.
Thus, we need to restrict the quantized flow \eqref{quant_euler} to $\SU(N)$, which is possible since $\SU(N)$ is a matrix Lie algebra (so it is closed under the matrix commutator $[\cdot,\cdot]$) and since the discrete Laplacian $\Delta_N$ restricts to an operator $\SU(N)\to\SU(N)$ (corresponding to the fact that the continuous Laplace--Beltrami operator $\Delta$ on $C^\infty(\Ss^2,\Cc)$ restricts to real functions $C^\infty(\Ss^2,\Rr)$).


Recall from the introduction that the continuous vorticity equation \eqref{euleq2} is a Lie--Poisson system with Hamiltonian given by \eqref{eq:ham_euler}.
Likewise, the quantized equation \eqref{quant_euler} is a Lie--Poisson system on the dual of $\mathfrak{su}(N)$ with Hamiltonian given by
\[
H(W)=\frac{1}{2}\Tr(\Delta^{-1}_N W W^\trans).
\]
The continuous Casimir functions $\mathcal C_k(\omega)$ for \eqref{euleq2} correspond, up to a normalization constant depending on $N$, to the following Casimir functions for \eqref{quant_euler}
\[
C_k(W)=\Tr(W^k) \mbox{  for } k=2,\ldots,N.
\]
As $N\to \infty$ we have convergence to the corresponding moments $\mathcal C_k(\omega)$ of the continuous vorticity \citep[see][Cor.~8.1.2]{RiSt2014}. 
We remark that the matrices $T^N_{lm}$, with the Frobenius inner product, share the orthogonality properties of $Y_{lm}$, with the $L^2(\Ss^2,\Cc)$ inner product. 
Therefore, if the initial vorticity $\omega$ is represented by a finite linear combination of spherical harmonics, then choosing $N$ sufficiently large, the continuous Hamiltonian $H(\omega)$ and enstrophy (quadratic Casimir) $\mathcal C_2(\omega)$ exactly coincide with the quantized analogs $H(W)$ and $C_2(W)$.


In the rotating case the quantized system is
\begin{equation}\label{eq:quant_euler_cor}
\dot{W}=[\Delta^{-1}_N (W-F),W]_N,
\end{equation}
where $F=2\Omega \mathrm{i} T^N_{10}$ represents the discrete Coriolis parameter.
The Hamiltonian in this case is given by
\[
H(W)=\frac{1}{2}\Tr(\Delta^{-1}_N (W-F) (W-F)^\trans).
\]
\subsection{Time discretisation} 

To obtain a complete algorithm we also have to discretize time.
For this, we use two different schemes. 
The first is implicit and preserves the Lie--Poisson structure.
The second is explicit but does not preserve the Lie--Poisson structure.

\subsubsection{Isospectral midpoint method (IsoMP)}
To take advantage of the quantization of the original equations, it is preferable to solve the quantized system \eqref{quant_euler} in time using a \emph{Lie--Poisson integrator}, i.e.\ a time-stepping scheme that preserves the Lie--Poisson structure \citep[cf.][]{McMoVe2014c,McMoVe2016b}. 
This way we obtain exact conservation of the Casimir functions and near conservation of the Hamiltonian (in the sense of backward error analysis of symplectic integrators \citep[cf.][]{HaLuWa2006}).
Since \eqref{quant_euler} is a Hamiltonian isospectral flow we can apply the Lie--Poisson schemes developed by \citet{modviv1}. 
We use here the second order isospectral midpoint rule (IsoMP). 
Given a time step parameter $h>0$ it is given by
\begin{equation}\label{IRK_mid}
\left. 
\begin{array}{ll}
&W_{n}=(I - \frac{h}{2}\Delta_N^{-1}\widetilde{W})\widetilde{W}(I + \frac{h}{2}\Delta_N^{-1}\widetilde{W})\\
&W_{n+1}=(I + \frac{h}{2}\Delta_N^{-1}\widetilde{W})\widetilde{W}(I - \frac{h}{2}\Delta_N^{-1}\widetilde{W}),
\end{array}
\right. 
\end{equation}
where $I$ is the identity matrix.
The matrix $\widetilde{W}$ is an auxiliary variable implicitly defined (together with $W_{n+1}$) by the two equations in \eqref{IRK_mid}. 
For further details on the method \eqref{IRK_mid} we refer to \citep{Viv2019}.

In presence of the Coriolis parameter $F$ the IsoMP scheme is
\begin{equation}\label{IRK_mid_Cor}
\left. 
\begin{array}{ll}
&W_{n}=(I - \frac{h}{2}\Delta_N^{-1}(\widetilde{W}-F))\widetilde{W}(I + \frac{h}{2}\Delta_N^{-1}(\widetilde{W}-F))\\
&W_{n+1}=(I + \frac{h}{2}\Delta_N^{-1}(\widetilde{W}-F))\widetilde{W}(I - \frac{h}{2}\Delta_N^{-1}(\widetilde{W}-F)).
\end{array}
\right. 
\end{equation}

The IsoMP method \eqref{IRK_mid} (and \eqref{IRK_mid_Cor}) exactly conserves angular momentum and the Casimirs $C_k(W)$, and nearly conserves the Hamiltonian $H(W)$ (its value oscillates in time without drift).

\subsubsection{Heun's method}

As an alternative to the Lie--Poisson preserving time discretization just described, we also consider the explicit Heun method.
Explicit methods, such as Heun's, exhibit linear drift in the first integrals.
However, if the linear drift is slow in comparison with the total simulation time, an explicit method might be the most competitive choice since it avoids non-linear root solving.
%
%
An efficient implementation of Heun's method for the quantized system \eqref{quant_euler} is the following:
\begin{equation}\label{Heun}
\left. \begin{array}{llll} 
K_1 = \Delta_N^{-1}W_n W_n  \vspace{.3cm}  \\
\widetilde{W} = W_n + h(K_1-K_1^\trans - \frac{1}{N}\Tr(K_1-K_1^\trans)I)\vspace{.3cm}  \\
K_2 = K_1 + \Delta_N^{-1}\widetilde{W}\widetilde{W} \vspace{.3cm}  \\
W_{n+1} = W_n + \frac{h}{2}(K_2-K_2^\trans - \frac{1}{N}\Tr(K_2-K_2^\trans)I).
\end{array}\right.
\end{equation}
In presence of the Coriolis parameter $F$ the scheme becomes
\begin{equation}\label{Heun_Cor}
\left. \begin{array}{llll} 
K_1 = \Delta_N^{-1}(W_n-F) W_n  \vspace{.3cm}  \\
\widetilde{W} = W_n + h(K_1-K_1^\trans - \frac{1}{N}\Tr(K_1-K_1^\trans)I)\vspace{.3cm}  \\
K_2 = K_1 + \Delta_N^{-1}(\widetilde{W}-F)\widetilde{W} \vspace{.3cm}  \\
W_{n+1} = W_n + \frac{h}{2}(K_2-K_2^\trans - \frac{1}{N}\Tr(K_2-K_2^\trans)I) .
\end{array}\right.
\end{equation}

\subsection{Complexity}\label{sub:complexity}

At first sight, it looks like the most demanding computational operation is the inversion of the discrete Laplacian $\Delta_N$:
it is a linear operator on $\SL(N,\Cc)$ and thus a fourth order tensor, so dense linear algebra would require $O(N^6)$ operations.
This is clearly not possible, even for moderate values of $N$.
However, from the formula \eqref{eq:discrete_laplacian} of Hoppe and Yau one can deduce
\begin{align*}
(\Delta_N)^{M_1'M_2'}_{M_1M_2} = &2\delta^{M_1'}_{M_1}\delta^{M_2'}_{M_2}
(s(s + 1) - M_1M_2)\\
&-\delta^{M_1'}_{M_1+1}\delta^{M_2'}_{M_2+1}\sqrt{s(s + 1) - M_1(M_1 + 1)}\sqrt{s(s + 1) - M_2(M_2 + 1)}\\
&-\delta^{M_1'}_{M_1-1}\delta^{M_2'}_{M_2-1}\sqrt{s(s + 1) - M_1(M_1 - 1)}\sqrt{s(s + 1) - M_2(M_2 - 1)},
\end{align*}
for $M_1,M_1',M_2,M_2'=1,\ldots,N$ and $s=(N-1)/2$. 
Notice that the tensor $\Delta_N$ is tridiagonal over the diagonal $M_1 = M_1'$ and $M_2 = M_2'$, i.e.\ it is sparse and contains only $O(N^2)$ non-zero entries; we store $\Delta_N$ as an $N^2\times N^2$ sparse matrix.
Remarkably, this sparse matrix also admits a sparse $LU$-factorization, i.e.\ a factorization of upper and lower diagonal matrices $L$ and $U$ which are also sparse with $O(N^2)$ non-zero entries.
Thus, to compute the inverse $\Delta_N^{-1}W$ requires just a single $LU$-factorization (which is $O(N^3)$ operations) and thereafter only $O(N^2)$ operations every time $\Delta_N$ is applied.
In essence, since the number of time steps for long-time simulations typically are of the order $O(10^6)$, this means that inversion of the discrete Laplacian only counts as $O(N^2)$ operations.



We solve the non-linear equation \eqref{IRK_mid} with Newton iterations. 
Thus, under the assumption that the average number of iterations per step is independent of $N$, the global complexity of the algorithm per time step is first $O(N^2)$ (for applying $\Delta_N^{-1}$) and then $O(N^3)$ (for the two matrix multiplications corresponding to computing the commutator $[\cdot,\cdot]$).
In summary, this means that the full complexity of the algorithm, per time step, is $O(N^3)$.



\subsection{Time scaling}\label{sub:real_time}
Recall that the correspondence between the matrix commutator on $\SU(N)$ and the Poisson bracket on $C^\infty(\Ss^2,\Rr)$ is $N^{3/2}[\cdot,\cdot] \approx c\lbrace\cdot,\cdot\rbrace$ for some constant $c$.
The requirement that 1 time unit of the vorticity equation \eqref{euleq2} correspond to 1 time unit of the quantized system \eqref{quant_euler} as $N\to \infty$ implies $c = \sqrt{16\pi}$.
In our simulations below we normalize the time scaling of the quantized equations by rescaling the initial conditions by $\lVert W_0\rVert$ and setting $[\cdot,\cdot]_N = [\cdot,\cdot]$.
This way, the non-dimensional time step $h$ correspond to $$\delta t = \dfrac{h\sqrt{16\pi}}{N^{3/2}\lVert W_0\rVert}$$ seconds of real time.
In all our simulations below we use the non-dimensional time step $h=0.1$.
A summary of the complete algorithm is given in Algorithm~\ref{alg:main}; it is implemented using MATLAB and available online.\footnote{The code is available at \href{https://bitbucket.org/kmodin/euler-sphere-quantization}{bitbucket.org/kmodin/euler-sphere-quantization}}




\begin{algorithm}
\caption{Summary of the complete numerical algorithm}\label{alg:main}
\noindent
\begin{algorithmic}[1]
    \Statex \textsc{Input}
    \State $\omega_0^{lm}$ \Comment{spherical harmonics components of the initial vorticity $\omega_0$}
    \State $N$ \Comment{quantization parameter}
    \State $h$ \Comment{non-dimensional time step parameter}
    \State $x_1,\ldots,x_{k}$ (optional) \Comment{nodes on the sphere for output}

	\vspace{2ex}  
    \Statex \textsc{Pre-processing}
    \State compute basis $T^N_{lm}$ to define the initial value $W_0\in \SU(N)$ from $\omega_{0}^{lm}$
    \State compute sparse $LU$-decomposition of the discrete Laplacian $\Delta_N$

	\vspace{2ex}         
    \Statex \textsc{Time-stepping}
	\For{$n\gets 1, n_{\it max}$} 
		   \State apply $\Delta_N^{-1}$ (using sparse $LU$-factorization) \Comment{$O(N^2)$ operations}
		   \State compute matrix commutator \Comment{$O(N^3)$ operations}
	\EndFor

    \vspace{2ex}  
    \Statex \textsc{Post-processing}
    \State extract the spherical harmonics components from the output states $W_n$
    \State use the nonequidistant spherical FFT transform~\citep{KePo2008} to recover $\omega(x_1),\ldots,\omega(x_k)$
\end{algorithmic}
\end{algorithm}

\section{Simulation results}\label{sec:simulations}
\subsection{Initial data with zero momentum}\label{sec:DQMsimulation}
We run our method with the same (randomly generated but zero momentum) initial data suggested by DQM, i.e.\ \citet{dri}.
We use $N=501$, $[\cdot,\cdot]_N=[\cdot,\cdot]$, and a dimensionless time step of $h=0.1$. 
With these parameters, the simulation time $t_k$ at step $k$ in the original units of time is computed by the formula $t_k = k/13643$,
(derived from the formula in \autoref{sub:real_time}).
We simulate with both the IsoMP and the Heun time integration.
For IsoMP, we use Newton-type iterations with a tolerance of $10^{-13}$.

\begin{figure}
\centering
IsoMP time integration \\
\includegraphics[width=0.99\textwidth]{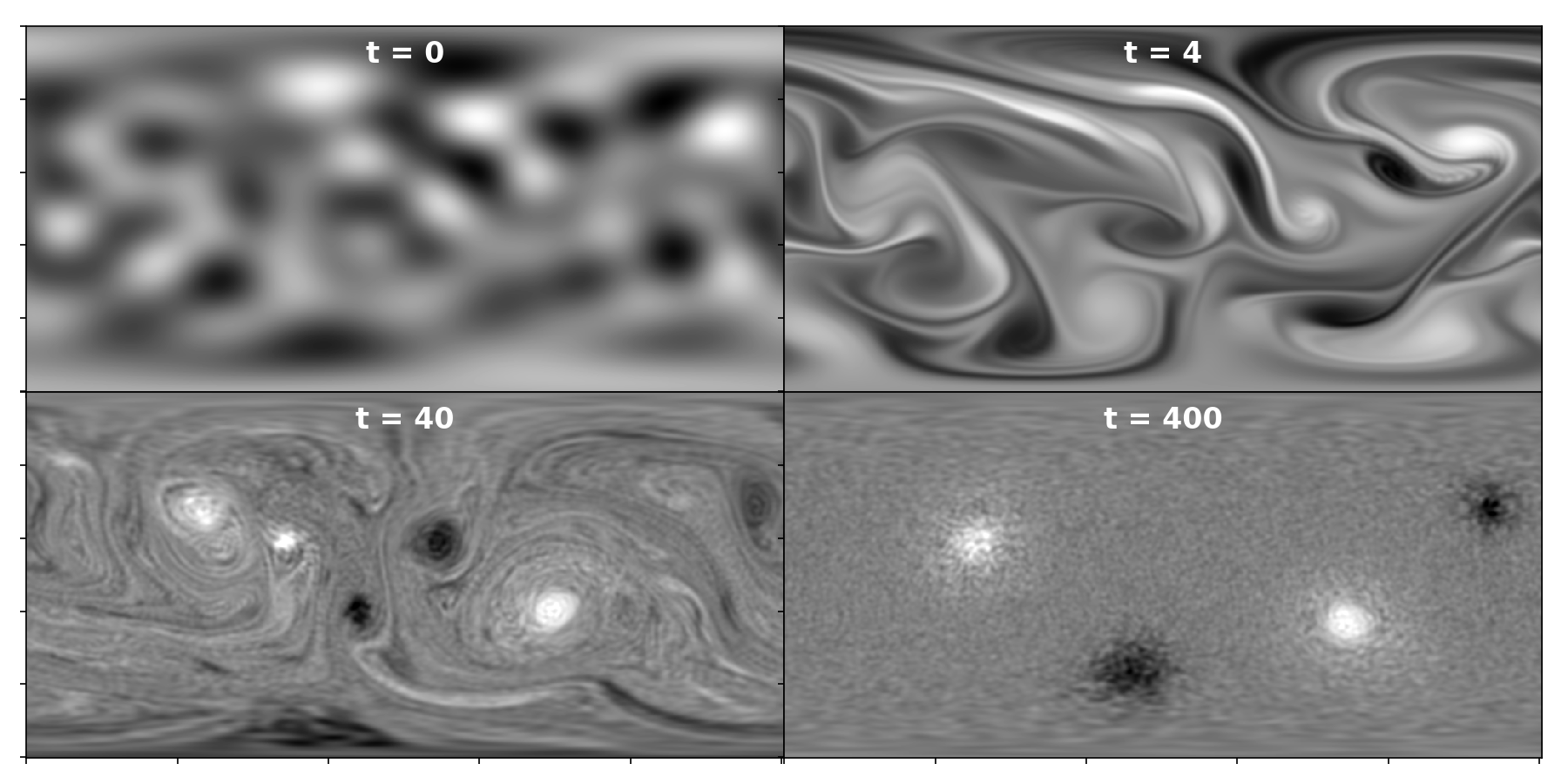} \\
Heun time integration \\
\includegraphics[width=0.99\textwidth]{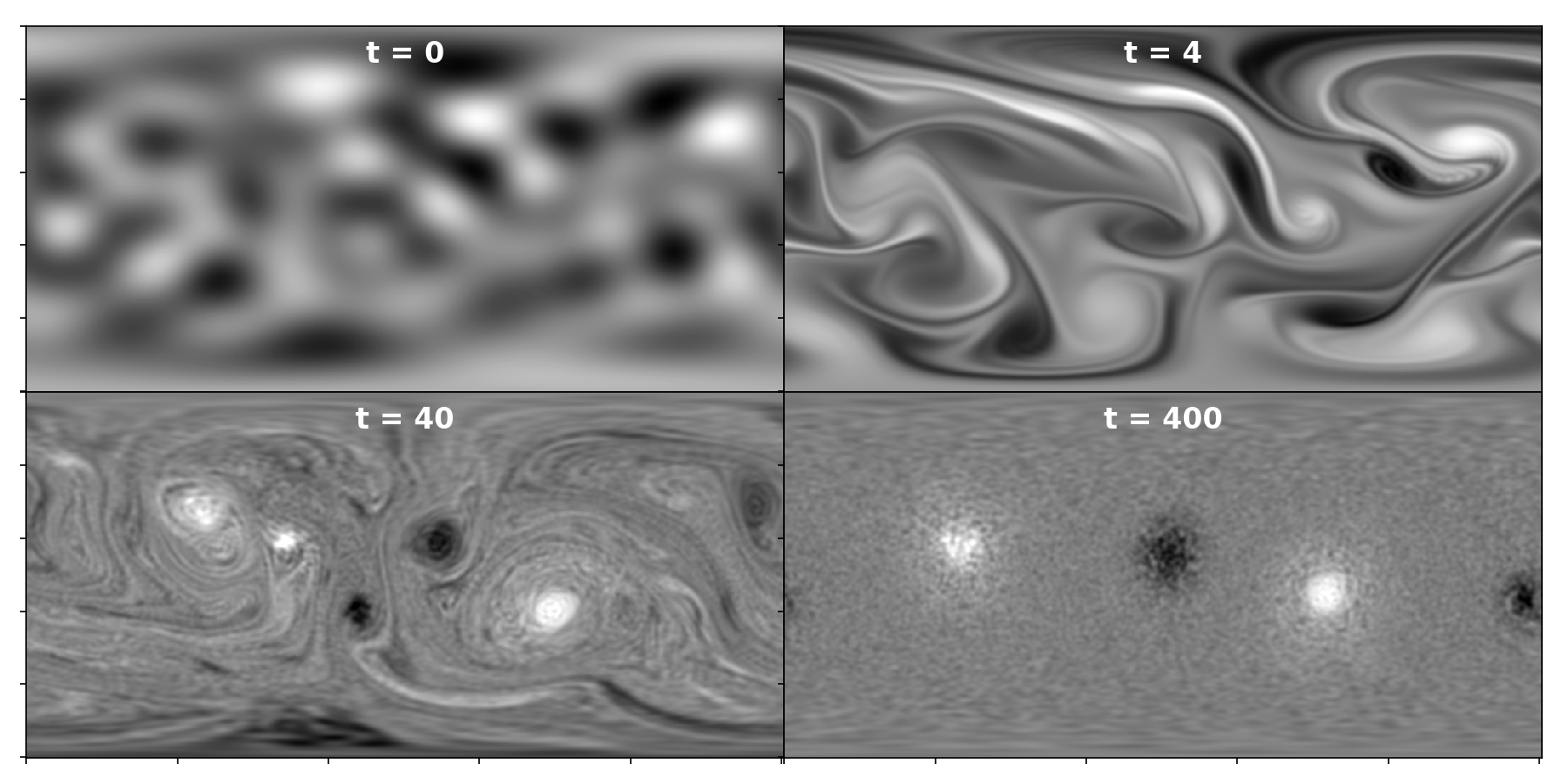}
 \caption{Simulation with two different time integration methods: IsoMP \eqref{IRK_mid} and Heun \eqref{Heun}. 
 Vorticity $\omega(x,t)$ clockwise from the top-left at $t=0s,4s,40s,400s$, for the initial data in DQM. The horizontal axis is the azimuth $\varphi\in[0,2\pi]$ and the vertical axis is minus the inclination $\theta\in[0,\pi]$.
 The results are visually indistinguishable up to $t=40$.
 At $t=400$ there are some differences in the positions of the vortex blobs.
 See also Movie 1, Movie 2, and Movie 3 of the supplementary material.
 }\label{fig:sim_1_isoMP}
\end{figure}

As already discussed in the introduction, the numerical results by DQM show that steady state is not reached, but rather four main vortex formations that move around the sphere, surrounded by smaller-scale vortices. 
Let us now compare with our results.
The vorticity at various output times is displayed, using spherical coordinates, in \autoref{fig:sim_1_isoMP} for the two different time integrations methods (IsoMP and Heun).


At time $t=4$ our simulations and those in DQM give visually indistinguishable results. 
At the early-intermediate vorticity, at time $t=40$, there is already a clear visible difference to DQM.
However, there is no visible difference between our two numerical time-integration schemes.
This indicates that, for time step lengths in the selected range, the choice of discretization in space, rather than time, dominates the numerical errors.

At $t=400s$ all simulation show the same qualitative feature: four large vortices moving about in the domain.
The exact position of the vortices are different between all the simulations (also between IsoMP and Heun).
There are two positive and two negative vortex blobs.
The exact strengths vary slightly between the blobs (see \autoref{sec:pvd} for further discussion about the vortex strengths).

When we run the simulation, either IsoMP or Heun, for long-times a clear pattern emerge: the 4 vortex formations are moving quasi-periodically.
The initial vortex mixing phase, up until the four vortex blobs have been formed at about $t=200$, is captured in Movie 1 of the supplementary material.\footnote{All the movies are also available at \href{https://bitbucket.org/kmodin/euler-sphere-quantization}{bitbucket.org/kmodin/euler-sphere-quantization}}\
The fast-forward Movie 2 of the whole simulation shows a short emerging phase of vortex mixing followed by a stable but unsteady large-scale quasi-periodic interaction of the four vortices.
In \autoref{sec:pvd} we discuss in detail the relation to stability of quasi-periodic point vortex solutions.
Movie 3 shows a simulation with the same initial conditions, but at the higher spatial resolution $N=1001$.
The qualitative behaviour is the same, with four vortex blobs forming and then circulating about each other in a quasi-periodic fashion.
However, the distribution of vortex formation is different in the high resolution simulation, with the positive instead of the negative blobs closer to the poles.

\begin{figure}
\includegraphics[width=0.49\textwidth]{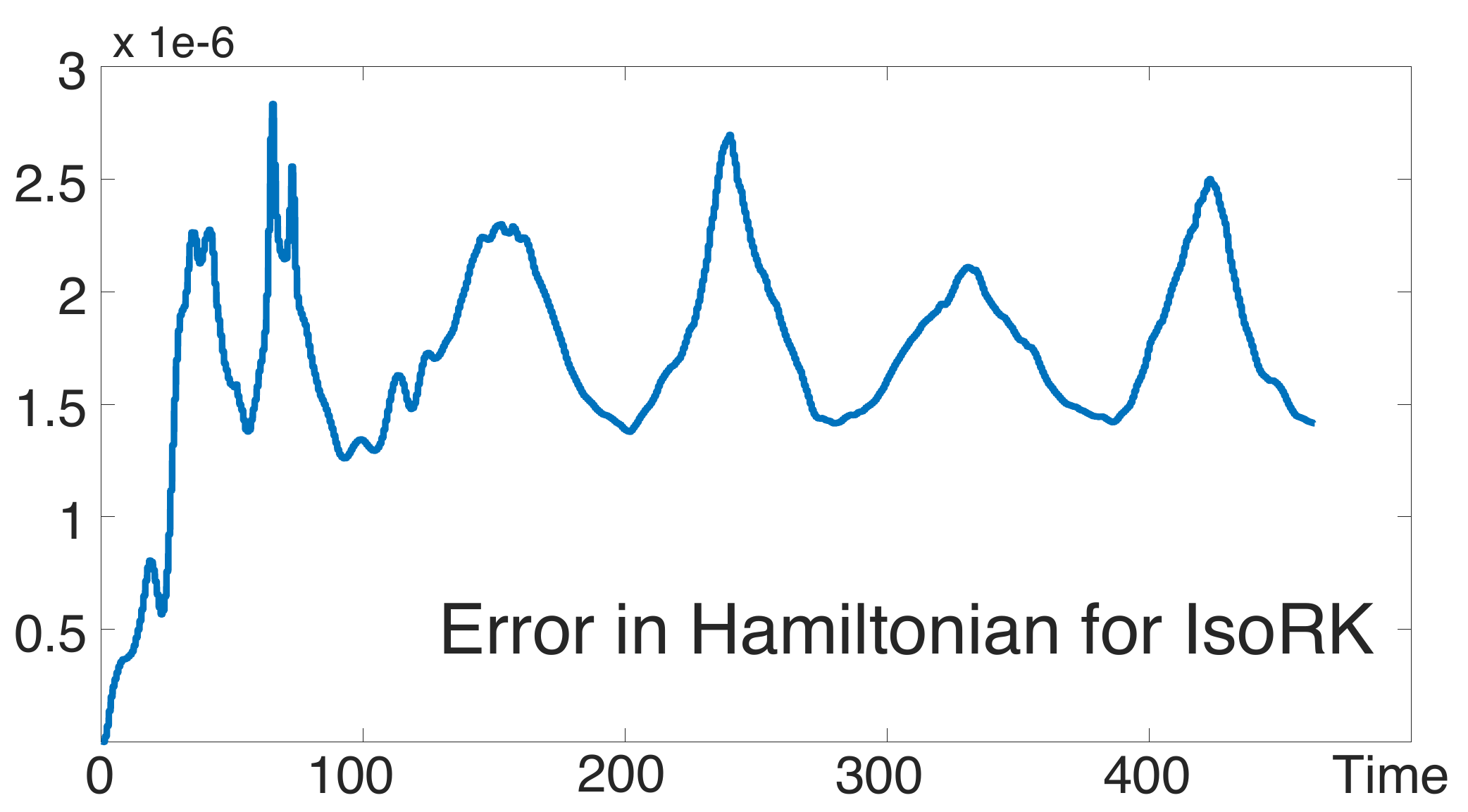}
\includegraphics[width=0.49\textwidth]{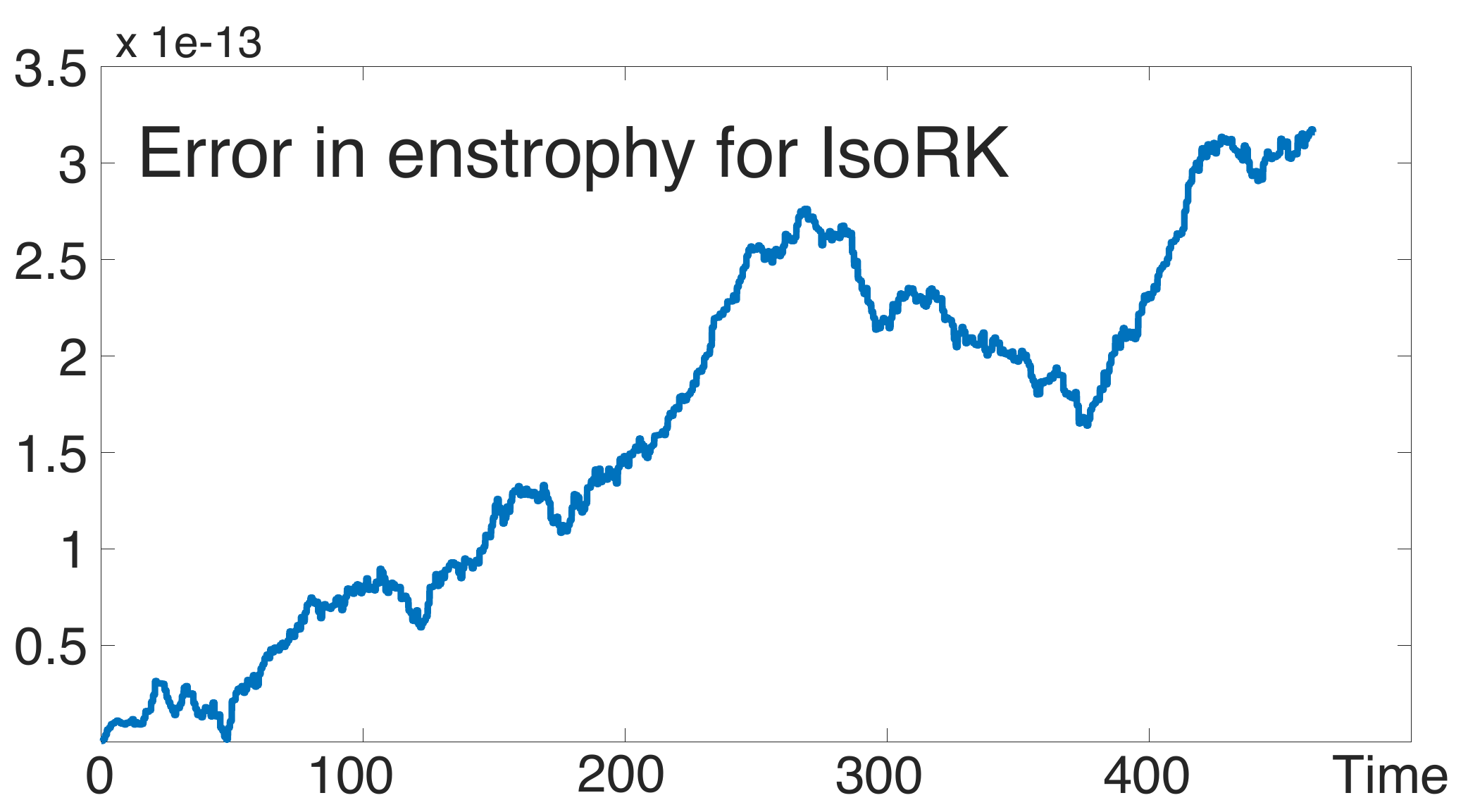}\\
\includegraphics[width=0.49\textwidth]{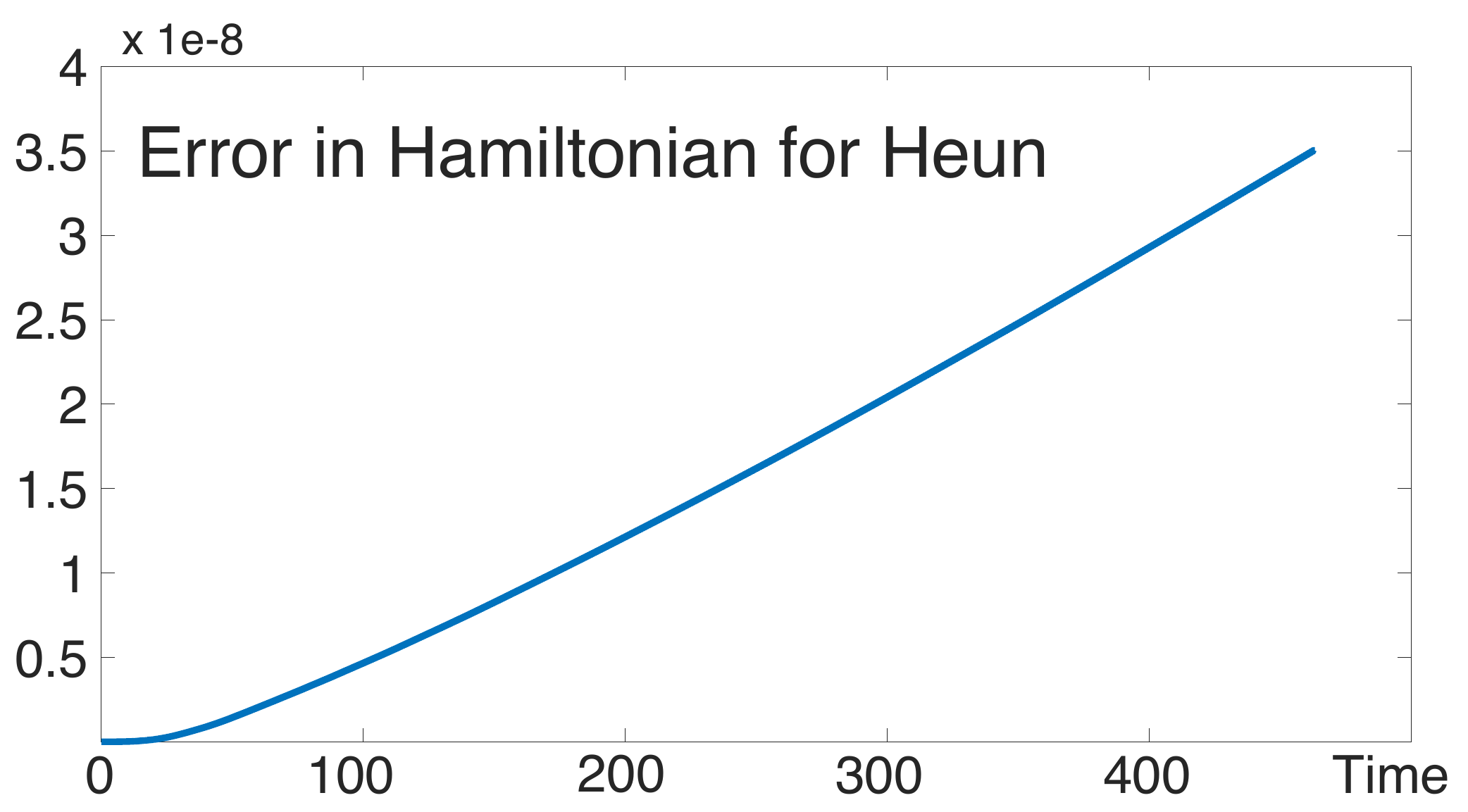}
\includegraphics[width=0.49\textwidth]{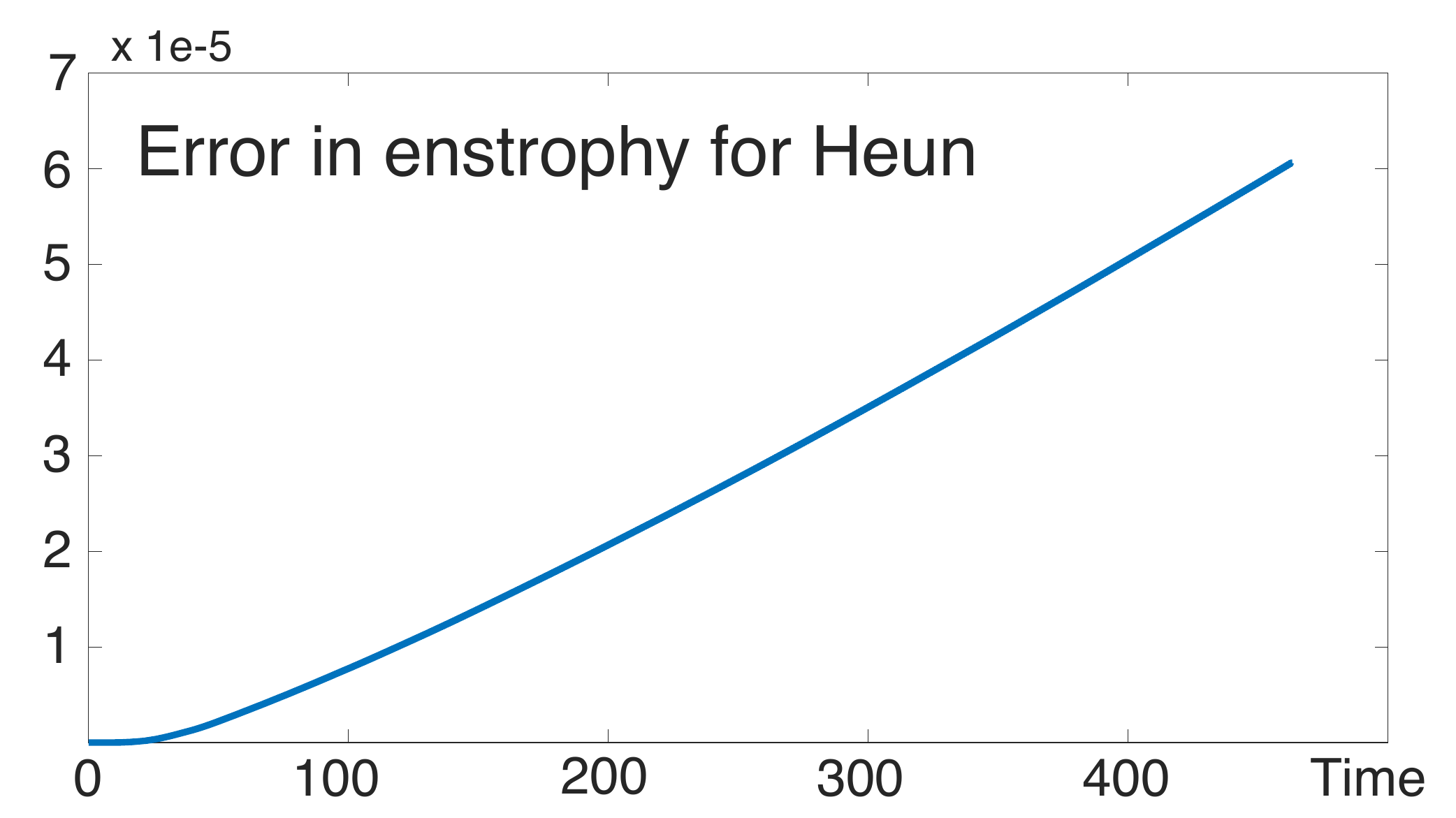} 
 \caption{Hamiltonian variation $|H-H_0|$ (left column) and enstrophy variation $|E-E_0|$ (right column) with the IsoMP (upper row) and Heun (lower row) time integrators.}
\label{fig:Ham_Ens_IsoMP}
\end{figure}

Let us continue the discussion here with the conservation properties of our method.
\autoref{fig:Ham_Ens_IsoMP} 
shows the variation of the energy and enstrophy during the simulation. 
For IsoMP, the energy is nearly conserved by a factor $10^{-6}$ with no sign of drift, whereas the enstrophy has the same variation as the Newton tolerance we have used, $10^{-13}$. 
For Heun, we see that albeit energy and enstrophy has a linear drift from their original values, the variation is quite small and in particular the energy changes less than with IsoMP. 
The negligible drift of energy and enstrophy is likely the reason why Heun perform so well.
We stress, however, that there is a drift, so at some point the numerics will break down, whereas with IsoMP such a breakdown will not occur since symplecticity is preserved.

The difference between IsoMP and Heun is more pronounced for the higher order Casimir functions of \eqref{quant_euler}. 
In fact, computing the maximal absolute variation of the eigenvalues of $W$, after $5\e{6}$ time steps, we get with IsoMP a value of the order $10^{-12}$, whereas with Heun a value of the order $1$. 
Even considering only the third and fourth momenta of the vorticity, the Heun scheme has an absolute variation, after $5\e{6}$ time steps, of the order $10^{-3}$.

In addition to integral invariants, such as energy and enstrophy, the continuous vorticity equation \eqref{euleq2} also conserves point-wise measures, such as the maximum vorticity $$\lVert\omega\rVert_{\infty} \coloneqq \sup_{x\in\Ss^2}\lvert\omega(x)\rvert.$$
Formally, the conservation of $\lVert\omega\rVert_{\infty}$ follows from conservation of the Casimir functions $\mathcal C_k(\omega)$ as $k\to \infty$.
Indeed, since the corresponding Casimir functions $C_k(W)$ of the quantized system approximate $\mathcal C_k(\omega)$ one can deduce (formally) that $\lVert\omega\rVert_{\infty}$ is nearly conserved without any drift (just like the energy).
In fact, this result follows rigorously from a theorem by \citet*[Thm.~4.1]{BoMeSc1994}, who proved that there is a constant $c\geq 0$, independent of $N$, such that 
$$
\lVert W\rVert \leq \lVert\omega\rVert_{\infty} \leq  \lVert W\rVert + \frac{c}{N}
$$ 
where $\lVert W\rVert$ is the matrix (operator) norm of $W\in \SU(N)$ and $\omega$ is the vorticity function corresponding to $W$.
Since $\lVert W\rVert$ is the largest eigenvalue (in magnitude) of $W$, and since all the eigenvalues are conserved by the quantized flow (the isospectral property), we get that $\lVert\omega\rVert_{\infty}$ is nearly conserved in the quantized system (i.e.\ it is an \emph{adiabatic invariant} for the quantized flow).




\begin{figure}
\centering
\begin{tikzpicture}
 \node (img) at (0,0)  {\includegraphics[width=0.5\textwidth]{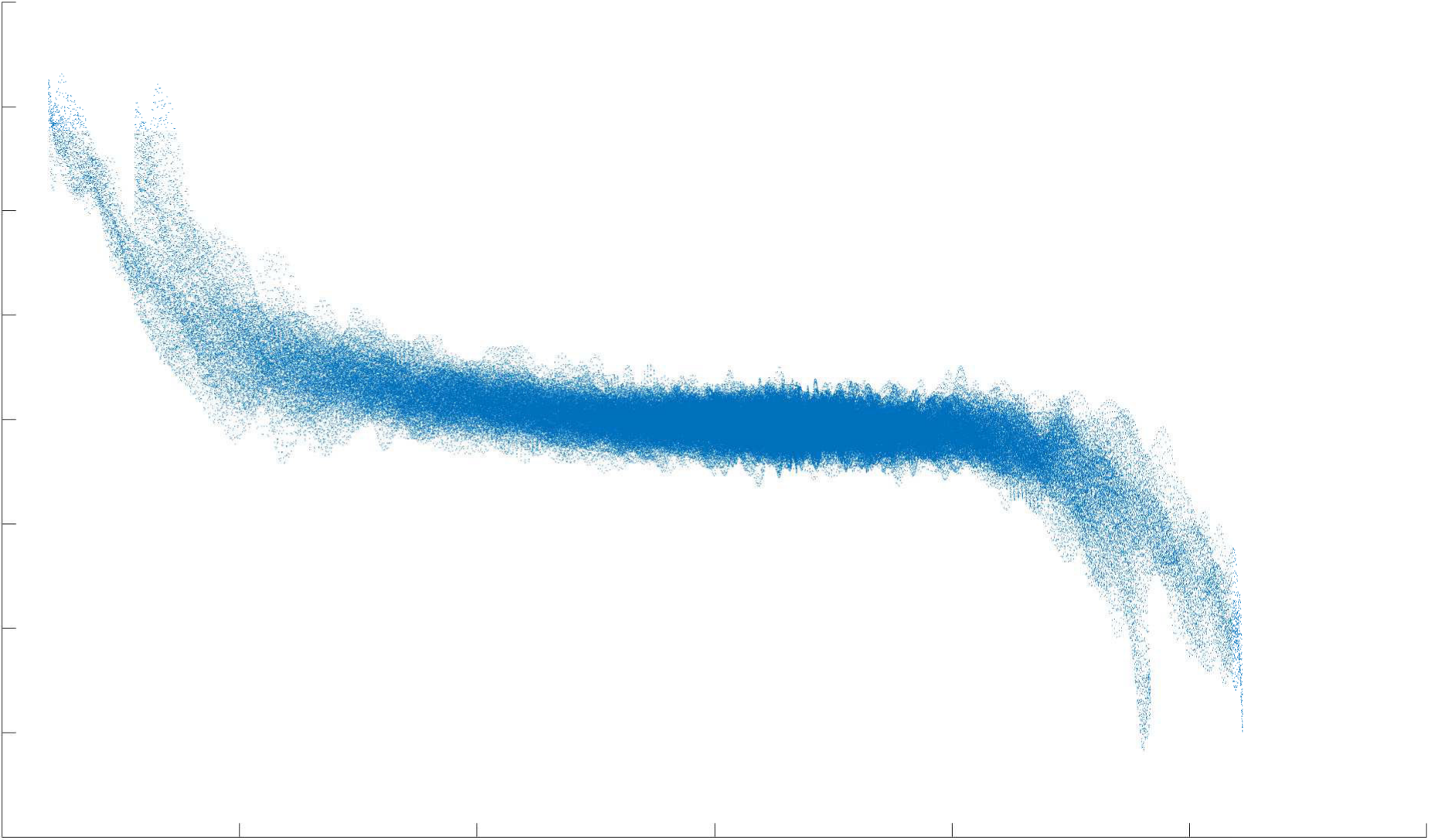}}; 
 \node (img2) at (0.5\textwidth,0) {\includegraphics[width=0.5\textwidth]{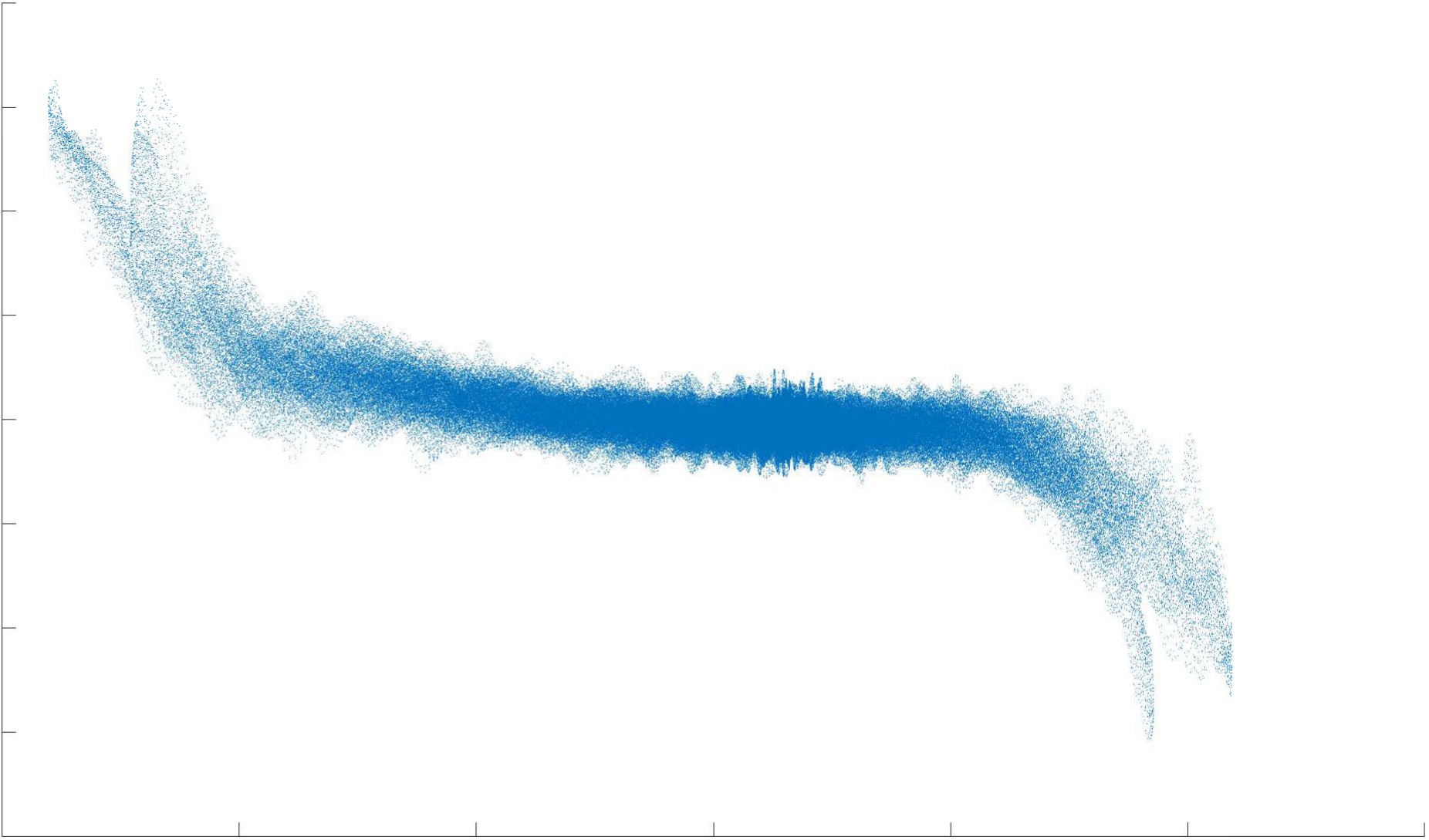}}; 
 \node at (-0.22\textwidth,0.0\textwidth) {$\omega$};
 \node at (0.0\textwidth,-0.12\textwidth) {$\psi$};
 \node at (0.28\textwidth,0.0\textwidth) {$\omega$};
 \node at (0.5\textwidth,-0.12\textwidth) {$\psi$};
 \node at (0,0.13\textwidth) {IsoMP time integration};
 \node at (0.5\textwidth,0.13\textwidth) {Heun time integration};
 \end{tikzpicture}
 \caption{Scatter plots of vorticity $\omega$ versus stream function $\psi$ at $t=400$. 
 (left) Using the IsoMP time integrator (\ref{IRK_mid}), and (right) using the Heun time integrator (\ref{Heun}).}
 \label{fig:W_PSI_400_Heun}
\end{figure}
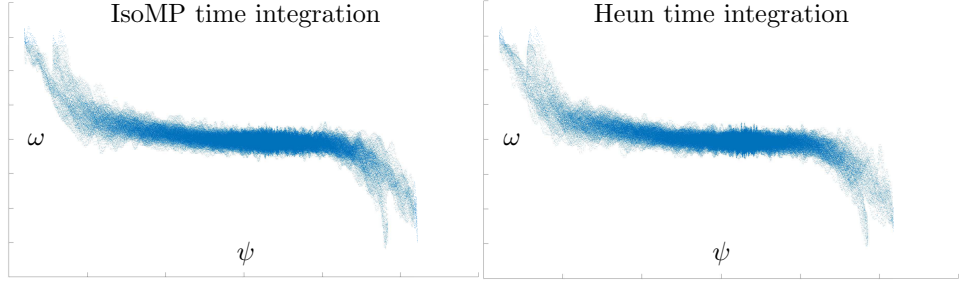

To measure the unsteadiness in the simulated flow we look at the relation between the vorticity $\omega$ and the stream function $\psi$ at $t=400$. 
The MRS theory predicts a steady flow determined by a functional relation $\omega = F(\psi)$ between the vorticity and stream function.
\autoref{fig:W_PSI_400_Heun} contains a scatter-plot of $\psi$ and $\omega$ for both IsoMP~\eqref{IRK_mid} and Heun~\eqref{Heun}. 
We notice that the shape of the resulting diagrams has branches, similar to those in DQM, indicating unsteadiness.
Our branches are more diffuse than those in DQM since no artificial dissipation is added in our model. 
We also see a slight difference between IsoMP and Heun: the one obtained with IsoMP has more defined branches. 

\subsection{Generic initial data}\label{sec:NonZeroMomsimulation}
In this section we present the results obtained with our numerical scheme on randomly generated initial conditions. 
We show that the generic behaviour for long times described by \citet{dri} it is not attained for non-zero angular momentum of the fluid. 
In our simulations we use $N=501$, $[\cdot,\cdot]_N=[\cdot,\cdot]$, and a dimensionless time step of $h=0.1$. 
Again, with these parameters the simulation time at step $k$ in the original units of time is computed by the formula $t = k/13643$.
We simulate with the the Heun time integration as it is faster;
for time evolutions as long as 400 real seconds the decay in enstrophy is negligible (see \autoref{fig:Ham_Ens_IsoMP}).

The generic random initial vorticity is obtained as it follows. 
Consider the expansion of the vorticity function in terms of spherical harmonics
\[
\omega(x) = \sum_{l=1}^\infty\sum_{m=-l}^l \omega^{lm}Y_{lm}(x).
\]
Then, $\omega\in L^2(\Ss^2)$ means that $\sum_{l=1}^\infty\sum_{m=-l}^l  | \omega^{lm}|^2<\infty$. 
We set the level of truncation $l_{\mbox{max}}=N-1=500$ and we generate the coefficients as random variables such that $\omega^{lm}l^{1 + \varepsilon}\sim\mathcal{N}(0,1)$, where $\mathcal{N}(0,1)$ is the normal Gaussian distribution and $\varepsilon = 10^{-3}$.
We stress that $L^2(\Ss^2)$ as the space for initial conditions is a natural choice in terms of Fourier analysis. 
Generating random initial conditions as just described corresponds mathematically to drawing samples from the isotropic Gaussian random field on $L^2(\Ss^2)$ as described by \citet{LaSc2015}.

\begin{figure}
\includegraphics[width=\textwidth]{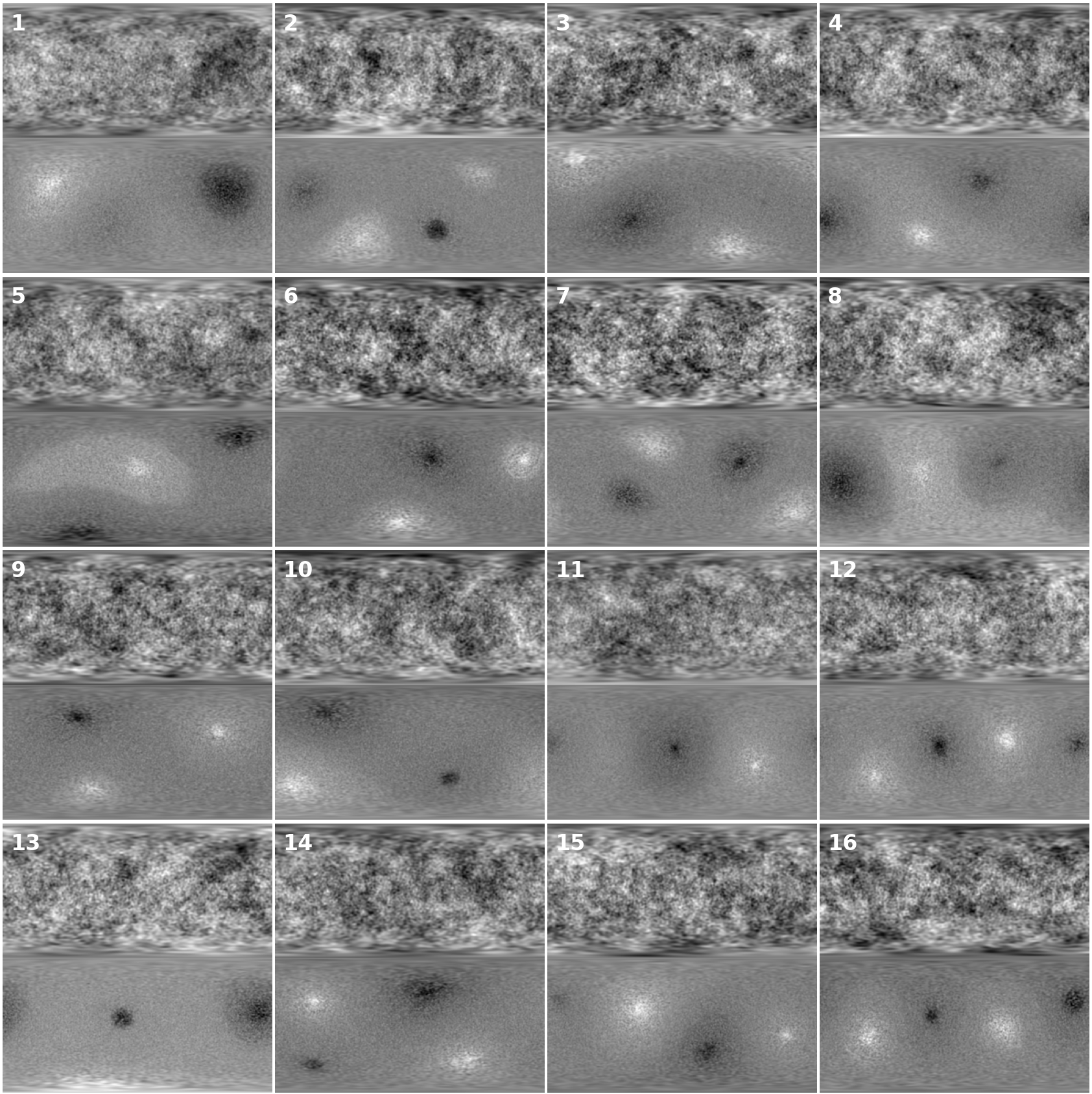}
 \caption{Pairs of initial (upper) and final (lower) vorticities for the 16 generic simulations with $L^2(\Ss^2)$ random initial data. 
 The numbers labelling the simulations correspond to those in \autoref{fig:ratio_momentum_enstrophy}.}
 \label{fig:l2_simulations}
\end{figure}

In this setting, we run 16 simulations on a cluster for long times. 
The vorticity for the simulations are given in \autoref{fig:l2_simulations} at time $t=0$ and $t=400$ real seconds. 
From the initially chaotic vorticity we see at time $t=400$ three qualitatively scenarios, with either 2, 3 or 4 persistent coherent vortices. 
In \autoref{sec:pvd} we explain this phenomenon in terms of integrability of point vortex dynamics and KAM theory. 
Movie 5, Movie 6, and Movie 7 of the supplementary material show the complete evolution of simulations 1 (giving 2 blobs), 4 (giving 3 blobs), and 7 (giving 4 blobs).

\section{Relation with point vortex dynamics}\label{sec:pvd}

We now explain the connection between the long-time behaviour of the Euler equations~\eqref{euleq1} on a non-rotating sphere and integrability theory of point vortex dynamics (PVD).
Recall from the introduction that our theory is based on the following two assumptions:
\begin{enumerate}
	\item The inverse energy cascade operates in such a way that smaller vortex formations of the same sign merge into larger ones when they get close enough.
	\item PVD describes the motion of vortex blobs well, as long as the blobs are well-separated so that no merging occurs.
\end{enumerate}
Based on the simplest, zero momentum case, we first give a detailed numerical verification of the second of these assumptions.
We then give the connection to integrability.
After that, we address the generic case of non-zero momentum and we show how our simulation results, with the three observed regimes, is a consequence of our theory. 

\subsection{Zero momentum case}\label{sec:pvd_zero_mom}
In this section we give a detailed study of the relation of our simulation results to the dynamics of four point vortices on the sphere, following up the brief study in~\citep{dri}.
For a detailed treatment of point vortex dynamics, we refer to the monograph of \citet{Ne2016} or the survey paper by \citet{Ha2007}.

Already \citet{Hel1858} knew that the incompressible Euler equations admit solutions with vorticity supported on single points.
Such solution also appear for the vorticity equations~\eqref{euleq2} on a sphere in the non-rotational case \citep{bog}.
That is, vorticity is a finite sum of $n$ Dirac delta distributions
\[
\omega = \sum_{i=1}^{n} \Gamma_i \delta_{x_i},
\]
where $\Gamma_i\in\Rr$ are the strengths and $x_i \in \Ss^2$ are the positions of the point vortices.
The solutions evolve according to an ordinary differential equations 
known as the \emph{point vortex equation}
\[
\dot{x_i} = \frac{1}{4\pi}\sum_{i\neq j} \Gamma_j\dfrac{x_j\times x_i}{1-x_i\cdot x_j},
\]
for $i=1,\dots,n$.
Notice that multiplying all $\Gamma_i$ by a factor does not change the phase space trajectories (only their speed), so only \emph{relative} strengths are of importance to us.
Our aim is to extract the positions and relative strengths of the vortex blobs in the DQM simulation, to compare their motion with the corresponding system of $n=4$ point vortices.

\begin{figure}
	\centering
	\includegraphics[width=0.99\textwidth]{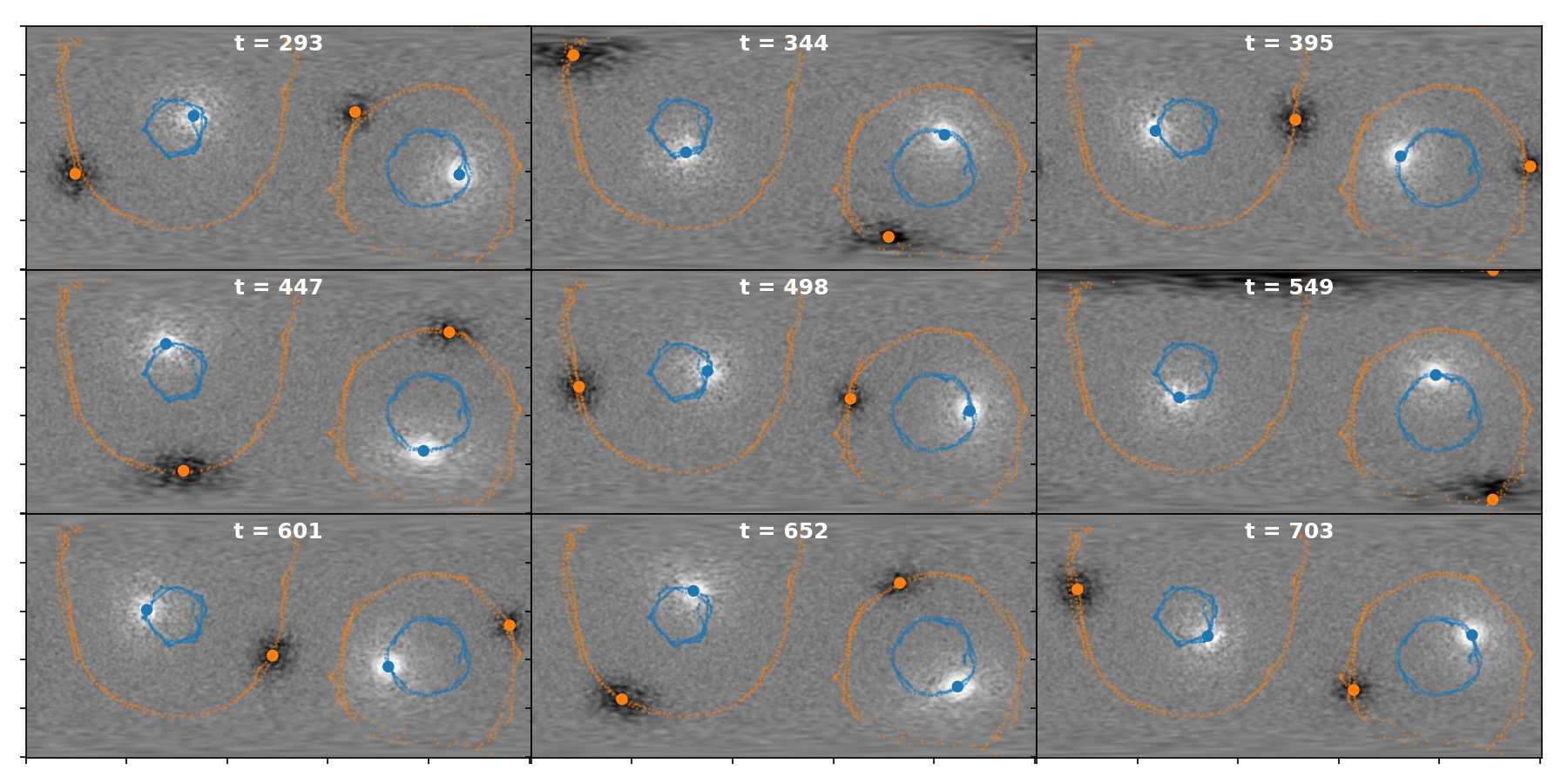}
	\caption{Behaviour of the vortex blobs for the same initial data as in \autoref{fig:sim_1_isoMP}.
	The quasi-periodic motion of the blobs is tracked in the scatter plots. 
	}
	\label{fig:DQMtracking}
\end{figure}

To extract the trajectories on the sphere of the 4 vortex blobs in the simulation from the previous \autoref{sec:DQMsimulation}, we use a tracking algorithm based on \texttt{Python/SciPy}.
The result is given in \autoref{fig:DQMtracking}.

Now, for the Euler equation on a non-rotating sphere, the total angular momentum $\mathbf{L}$ is conserved
\[
\frac{d}{dt}\mathbf{L} = \frac{d}{dt}\int \omega \mathbf{n} = 0.
\]
The DQM simulation is set up with vanishing total angular momentum, $\mathbf{L} = 0$.
Thus, the point vortex solutions should fulfill
\[
\int \sum_{i=1}^{n} \Gamma_i \delta_{x_i} \mathbf{n} = \sum_{i=1}^n \Gamma_i x_i = 0.
\]
If we set $\Gamma_1 = 1$ (since we are only looking for relative strengths), then, for generic positions $x_i$, this yields a linear set of equations from which $\Gamma_2,\Gamma_3,\Gamma_4$ can be determined from the positions alone.
The computed relative strengths thereby obtained correspond well with those obtained by numerical integration over circular domains covering the blobs.
In summary, we have the following extracted positions (expressed in inclination $\varphi$ and azimuth $\theta$) and corresponding computed relative strengths
\begin{equation}\label{eq:pvd_data}
\begin{aligned}
\varphi &= [2.3218,  -0.9638, -2.5283,  0.8511]\\
\theta &= [1.3017, 1.8837, 1.577, 1.5896] \\
\Gamma&=[1, 0.9002, -0.5436, -0.4178].
\end{aligned}	
\end{equation}

\begin{figure}
\includegraphics[width=0.99\textwidth]{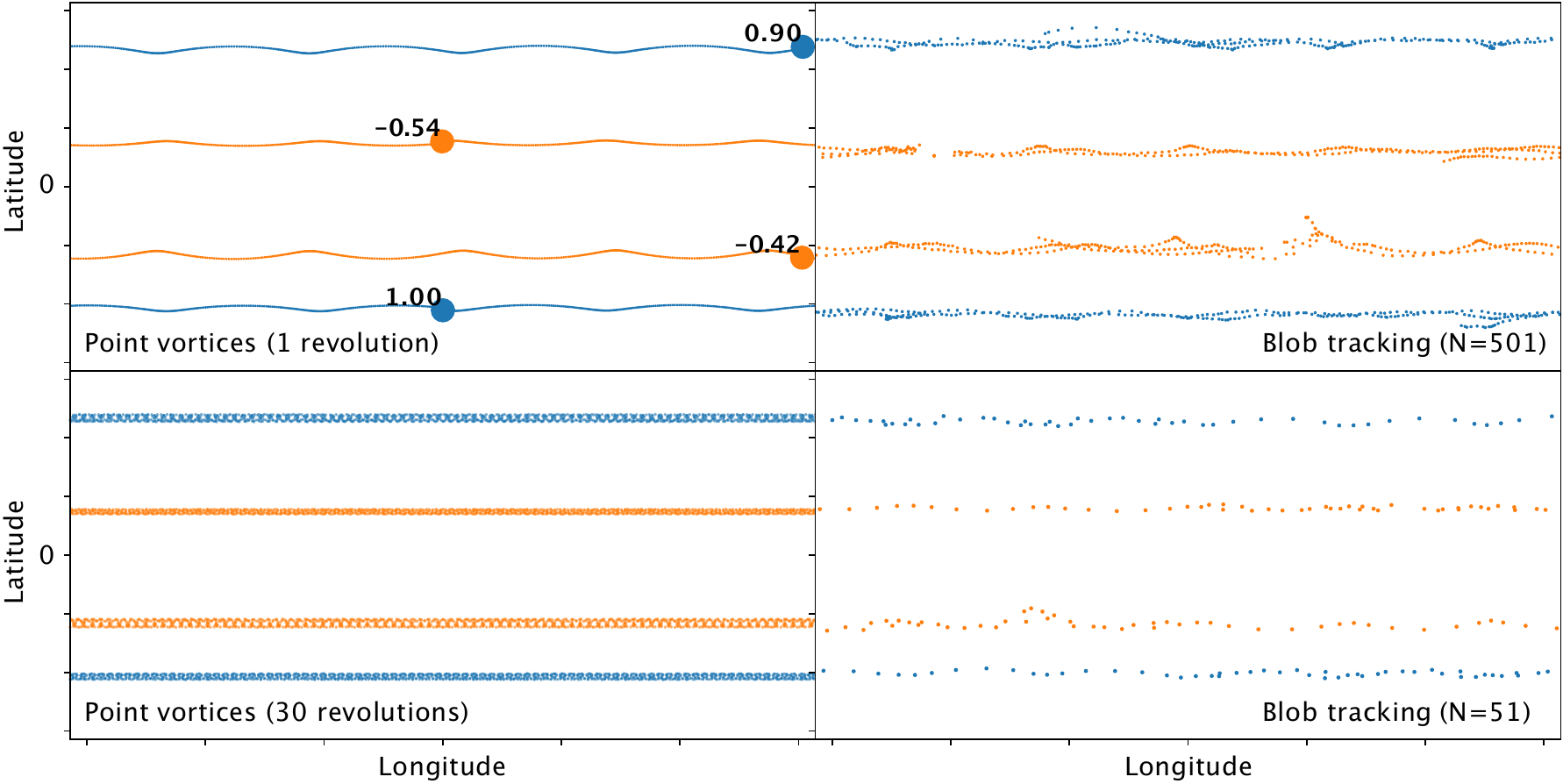} 
\caption{
Motion of point vortices (left) versus tracking of vortex blobs in the simulation in \autoref{fig:DQMtracking} (top-right) and in the simulation at low resolution (N=51), as described in (\ref{Gauss_blobs}) (bottom-right).
}
 \label{fig:blobs_PV}
\end{figure}

\begin{figure}
	\centering
	\includegraphics[width=0.5\textwidth]{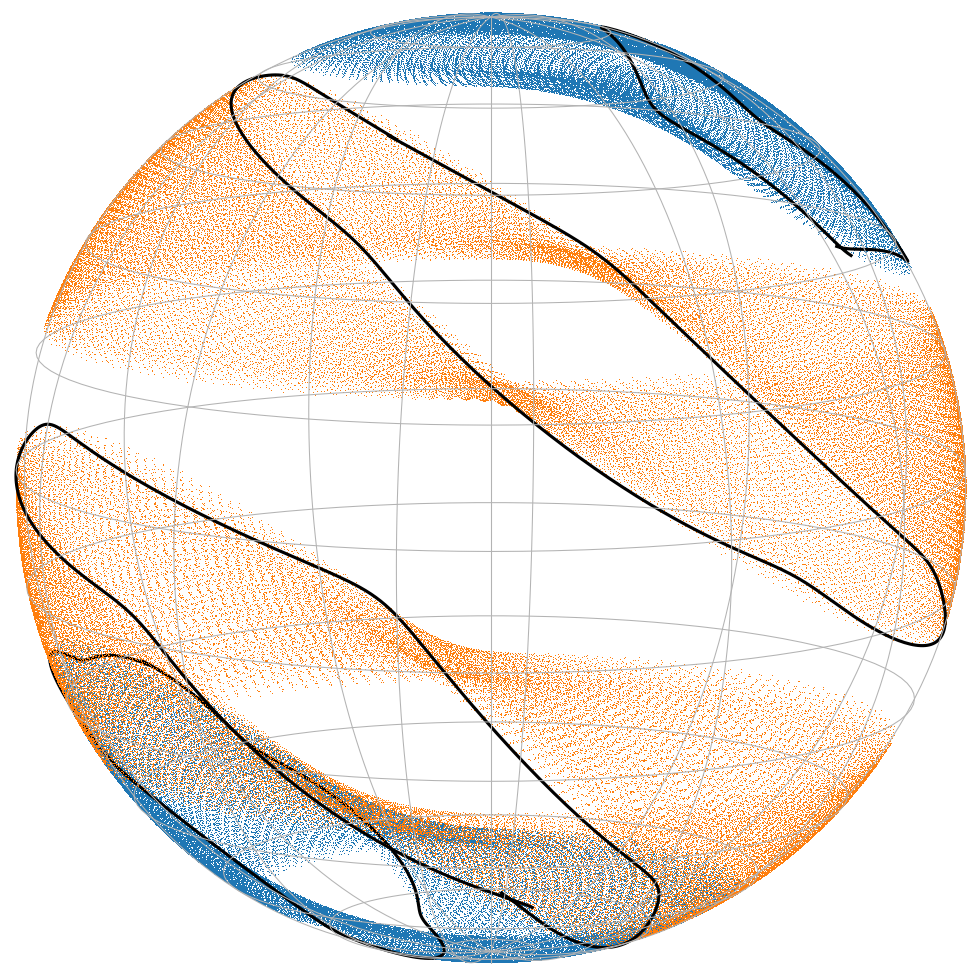}
	\caption{Long-time behaviour (about 2800 revolutions) of the four point vortex simulation. 
	The axis of rotation is slowly drifting, yielding a precession of the trajectories.
	The indicated lines show the motion in the last simulated period.
	These numerical results verify that the motion is integrable, i.e.\ quasi-integrable, with three or four frequencies: the wobbling during each rotation (highest frequency), the rotation about the axis (intermediate frequency), and the precession of the rotation axis (lowest frequency/frequencies).
	}
	\label{fig:pv_sphere_long}
\end{figure}

\begin{figure}
\includegraphics[width=0.99\textwidth]{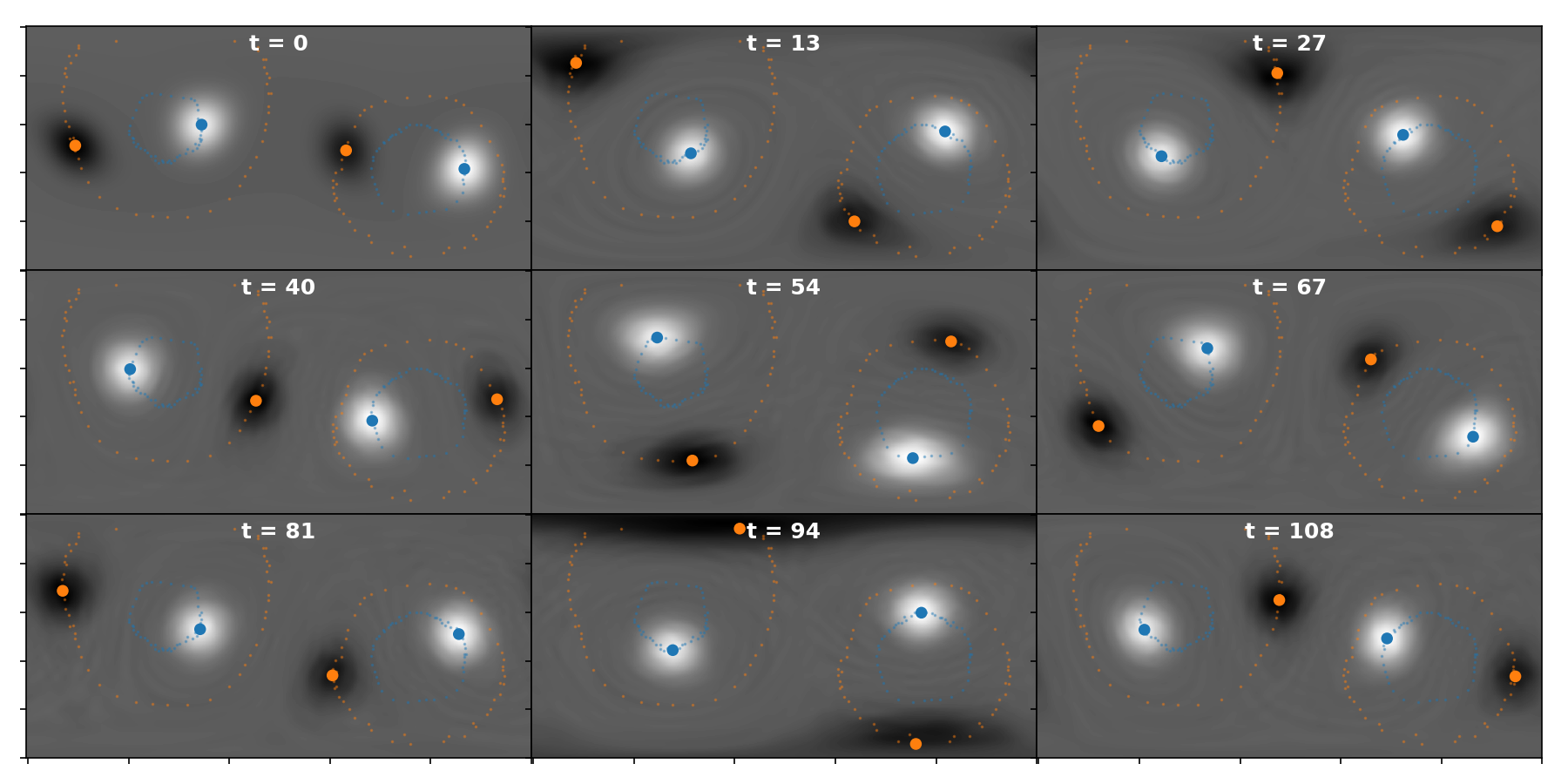}
\caption{
Motion of initially Gaussian blobs discretized with low spatial discretization ($N=51$).
The shape, strengths, initial positions of the blobs are given by \eqref{Gauss_blobs}.
The tracked motion of the blobs is in good agreement with the $N=501$ simulation and the point vortex dynamics (compare with \autoref{fig:blobs_PV} and the lower right diagram in \autoref{fig:blobs_PV}).
}
 \label{fig:4blobs}
\end{figure}
To obtain the absolute strengths, i.e., to determine the scaling, one might use the total energy integral, noting that the point vortex Hamiltonian is quadratic in the strengths.

\begin{rem}
The fact that the relative strengths of the point vortices are uniquely determined (in the generic case) by the positions given vanishing total angular momentum is interesting.
It shows that there is a connection between the strengths and the positions.
One may ask to what extent the strengths determine the positions for zero momentum configurations.
That is, what is the dimension of the manifold of four point vortices with vanishing total angular momentum.
Heuristically, just counting constraints, one gets the dimension of all possible four point vortex configurations, 8, minus the dimension of the 3 angular momentum constraints, which gives 5 dimensions.
To investigate this question in detail, one can use \emph{symplectic reduction theory} \citep[cf.][]{marat}.
\end{rem}

We run the point vortex dynamics simulation with data from \eqref{eq:pvd_data} using the symplectic Lie--Poisson integrator by \citet{McMoVe2016b}.
Let us now compare this point vortex simulation with the tracked blob motion in \autoref{fig:DQMtracking}.

In the chosen spherical coordinates, the motion of the tracked blobs in \autoref{fig:DQMtracking} looks complicated, but, in fact, when plotted on the sphere, one can see that it is almost a steady rigid rotation about a fixed axis.
By least squares we find the best approximating rotation axis, and we use new spherical coordinates with the new rotation axis as the north pole.
The resulting trajectories, of both the tracked blobs and the computed four point vortex dynamics with data~\eqref{eq:pvd_data}, are given in \autoref{fig:blobs_PV}.
We see that the motion between the point vortices and the tracked blobs are in good agreement, as also reported by \citet{dri}.

Looking at the almost pure rotational trajectories in \autoref{fig:blobs_PV}, it is natural to ask if there is an underlying relative equilibrium, i.e., a close-by solution given by a simultaneous, steady rotation of all the four point vortices.
The answer is no: any such relative equilibrium must in fact be a static equilibrium, because the total angular momentum is zero.
Thus, the `wobbling' in \autoref{fig:blobs_PV} is necessary for unsteadiness to occur.
Continuing this train of thought, we may look for static, non-stable equilibria for zero momentum four point vortex dynamics with arbitrary strengths.
Based on symmetry considerations, a general study of equilibria for point vortex dynamics on the sphere is carried out by \citet*{LaMoRo2011}.
For general strengths there does not seem to exist equilibria, but for pairs of two equal positive and two equal negative strengths there are, given by staggered rings \citep[see][Sec.~8]{LaMoRo2011}.
Since the computed strengths \eqref{eq:pvd_data} almost come in such pairs, and since at any instance in time the configuration of the vortex blobs in \autoref{fig:DQMtracking} is almost given by such staggered rings, we are, in this sense, always close to equilibria, but they are unstable.
That the strengths of the vortex blobs almost comes in pairs is likely not a coincidence.


The simulation in \autoref{sec:DQMsimulation} generating the blob formation is long enough to cover about 3 revolutions of the blobs about each other.
Although this is considered a `long-time' simulation of the Euler equations, it is not very long if one wants to study the stability of the quasi-periodic trajectories.
If we run the point vortex simulation for about 30 revolutions, we see in the lower left plot of \autoref{fig:blobs_PV} that the trajectories appear to keep on wobbling about the zonal lines.
But even 30 revolutions is not much.
With a much longer simulation of about 2800 revolutions, we see, as plotted in \autoref{fig:pv_sphere_long}, a different pattern emerge:
the positions of the vortices are spreading out by a very slow precession.
These numerical experiments indicate that the dynamics of the four point vortices restricted to the submanifold of vanishing total angular momentum is integrable, or at least quasi-integrable in the KAM sense. 
The frequencies would then be the oscillations within each revolution (highest), the rotation (intermediate), and one or two much lower frequencies for the precession.
In general, the dynamics of four point vortices is \emph{not} integrable.
However, the submanifold of vanishing total angular momentum is special, as it is the only submanifold of fixed angular momentum that is invariant under arbitrary rotations (if you rotate a configuration of zero momentum, it still has zero momentum). 
Indeed, \citet{Sa2007} showed that the dynamics of four point vortices with zero angular momentum is integrable.
As a theoretical approach aiming to prove integrability, one could also proceed by \emph{zero momentum Hamiltonian reduction} \citep[cf.][]{MaMiOrRa2007}.
Roughly, it goes as follows.
The Lie group $\mathrm{SO}(3)$ of rotations acts on the configuration space $(\Ss^2)^4$ of point vortices.
The corresponding Nöther integrals are the total angular momentum.
Since the area form on $\Ss^2$ is preserved by rotations, the action is symplectic.
By \emph{Poisson reduction} we thereby obtain a new Hamiltonian system on the Poisson manifold $(\Ss^2)^4/\mathrm{SO}(3)$ of dimension 5.
Now, every Poisson manifold is foliated in symplectic leafs.
In particular, we have the special \emph{zero momentum leaf}, given by restriction to the zero set of the total angular momentum.
We thereby obtain a Hamiltonian system on the symplectic manifold given by the zero momentum leaf of dimension 2, which is always integrable. 

Our findings in this section show that the initial conditions used by DQM, although random in the higher order spherical harmonics, is special since it has zero angular momentum.
That is, one cannot expect the long-time behaviour obtained with the DQM initial conditions to be generic for initial conditions with non-zero angular momentum.
Indeed, if four vortex blobs in a non-zero momentum configuration are formed, there might be further mixing, since their motion most likely will be chaotic.
For zero momentum, however, the quasi-periodic behaviour acts as a barrier, preventing further mixing.
We thus predict that for vanishing angular momentum, quasi-periodic asymptotics is the generic behaviour.
To investigate this question in detail is yet another future topic.




We now want to illustrate how the quasi-periodic motion of the blobs can be obtained even for a very coarse spatial discretization $N=51$.
The initial vorticity here is:
\begin{equation}\label{Gauss_blobs}
\omega_0(x) = \sum_{i=1}^4 \Gamma_i\exp(-20|x-x_i(\varphi_i,\theta_i)|^2)+C(x),
\end{equation}
for $\Gamma,\varphi,\theta$ as in \eqref{eq:pvd_data}, and $C(x)$ such that $\omega_0$ integrate to zero and has momentum $\mathbf{L}=0$.
The result is given in \autoref{fig:4blobs} and Movie 4, and the resulting vortex blob motion tracking (in adapted spherical coordinates) is given in the lower right plot in \autoref{fig:blobs_PV}.
We obtain good agreement with both the point vortex simulation and the full high resolution simulation with $N=501$.
That Gaussian vortex blob simulations can be carried with small discretization parameters $N$ is important, because it opens up for much longer simulations studying the stability of quadruple vortex blob formations.
We anticipate that the slow precession seen in point vortex dynamics also happen in vortex blob dynamics.

\subsection{Generic case}\label{sec:pvd_non_zero_mom}
\begin{figure}
	\centering
	\includegraphics[width=1\textwidth]{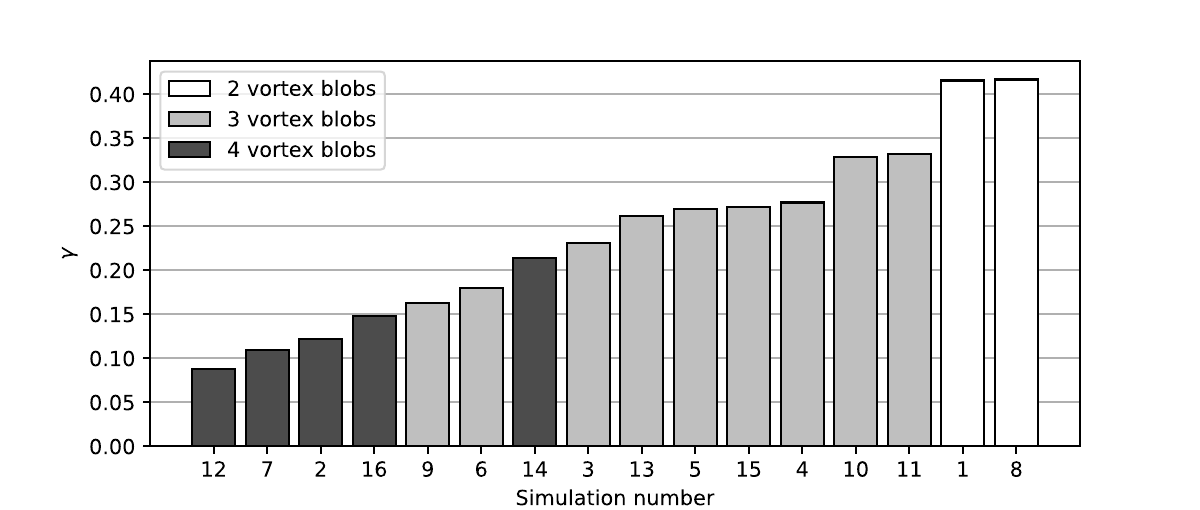}
	\caption{Values of $\gamma=\| \mathbf{L}\|/(R\sqrt{\mathcal{C}_2})$ for the simulations of \autoref{sec:NonZeroMomsimulation}. 
	The grey-scale correspond to the number of vortex blobs observed in the final state: 2, 3, or 4.
	Notice that the value of $\gamma$ largely determines the number of vortex blobs in the final state: 4 when $\gamma\lesssim 0.15$, 3 when $0.15\lesssim\gamma\lesssim 0.4$, and 2 when $\gamma\gtrsim 0.4$.
	} 
	\label{fig:ratio_momentum_enstrophy}
\end{figure}

The ideas presented in the previous paragraph can be extended to the non-zero momentum vorticity.
Our simulations in \autoref{sec:NonZeroMomsimulation} suggest that the four blobs formation in \citep{dri} is specific for initial conditions with small angular momentum. 
In fact, as already mentioned in the introduction, there is a correlation between the first integral $\gamma:=\| \mathbf{L}\|/(R\sqrt{\mathcal{C}_2})$ and the number of coherent vortices that persist in the final state (see \autoref{fig:ratio_momentum_enstrophy}).
As can be seen in the simulations, there exist a finite range for $\gamma$ (approximately $0.15\lesssim\gamma\lesssim 0.4$) for which the mixing of vortex blobs continuous until the dynamics reach a quasi-periodic motion of three blobs and no more mixing occurs after that.
Above this range, the momentum prevails on the other modes, allowing only the persistence of two large vortices. 
To the best of our knowledge this phenomenon has not been previously described. 

Based on the two assumptions presented at the beginning of this section, an explanation for the observed phenomenon is offered through integrability properties of point vortex dynamics as already laid out in the introduction.
Indeed, it is known that for non-zero momentum the three point vortex dynamics is integrable, whereas it is not integrable in general for four point vortices~\citep{Sa2007}. 
In \autoref{sec:NonZeroMomsimulation} our numerical simulations show that when the angular momentum $\mathbf{L}$ is non-zero there may occur further mixing from the four vortices found by \citet{dri}, leading to a final state of three or two vortices. 
This can be understood in terms of perturbation of an integrable configuration of point vortices. 
In fact, starting from the zero momentum four blobs, one can understand the modification of momentum $\mathbf{L}$ to a non-zero value as a perturbation of the zero-level set. 
As noticed in \autoref{sec:NonZeroMomsimulation}, up to the critical value of $\gamma\sim 0.15$, the four point vortex dynamics persists and is quasi-periodic.
The reason for such a situation is that the momentum $\mathbf{L}$ is a small perturbation, in the sense of the KAM theory \citep[cf.][]{Se1995}, of an integrable system of point vortices, and small perturbations do not destroy the invariant tori, so the quasi-periodicity is still intact, acting as a barrier for further mixing.
However, when $\gamma\gtrsim 0.15$, the perturbation from zero momentum is large enough to destroy the four point vortex integrable state, leading to chaotic trajectories of the blobs and therefore further mixing up to the next integrable configuration of three point vortices. 
Eventually, increasing the magnitude of the momentum $\mathbf{L}$ over a certain threshold ($\gamma\gtrsim 0.4$), the final state of the vorticity can be described by two antipodal point vortices only, aligned in the direction of the momentum $\mathbf{L}$.
These two vortices are now so large from the start that they tend to directly swallow the smaller vortices without passing through the quasi-periodic three blobs formation (although, if one looks carefully in simulations 1 and 8 corresponding to the higher values of $\gamma$, one can trace a small third vortex blob which does not affect the dynamics).

\begin{rem}\label{rem:torus}
We point out that the relation of the final state of the total vorticity and the point vortex dynamics strongly depends on the manifold where the equations take place. 
On a torus, in fact, the total angular momentum does not play any role for the integrability of the point vortex dynamics. 
Instead, the total circulation of the point vortices (i.e. the sum of the vortices' strengths) is determinant. 
With zero circulation, three point vortex dynamics on the torus is integrable, whereas for non-zero circulation only two point vortex dynamics is integrable. 
This explains why the latter configuration of two large blobs has been extensively observed \citep[see for example][]{QiMa2014}, whereas no three large blobs on a torus appear in the simulations. 
Indeed, prescribing a final state of zero circulation of three blobs (notice, not only a zero circulation of the total vorticity since one has to subtract the constant background circulation) is not possible, unlike prescribing zero momentum. 
Hence, on a torus, our theory predicts that the generic behaviour is the formation of two steady point vortices, as also observed numerically.
We hypothesize that our theory, connecting long-time behaviour with integrability of point vortex dynamics, is valid for the Euler equations on any 2D surface.
To investigate, and possibly verify, this claim in full detail is a future research topic.
\end{rem}

\section{Conclusions and outlook}\label{sec:conclusions}
We have developed a new numerical algorithm for the Euler equations on a sphere that preserves, up to machine precision, the Casimir functions of \eqref{quant_euler}-\eqref{eq:quant_euler_cor} and nearly conserves the Hamiltonian (see section~\ref{sec:num_scheme}).
The spatial discretization is based the work of \citet{hopPhD} on geometric quantization of the infinite dimensional Poisson algebra of smooth functions on the sphere by matrix Lie algebras $\SU(N)$ for an increasing~$N$ (corresponding to the spatial discretization parameter), and the suggestion of \citet{ze2} to use this quantization for Euler fluids.
The resulting finite dimensional dynamical system on $\SU(N)$ is isospectral, corresponding to conservation of Casimir functions, and preserves a Lie--Poisson structure, corresponding to the Hamiltonian structure of ideal fluids.

On the one hand, long-time simulations on a non-rotating sphere for zero momentum initial vorticity confirm the results in \citep{dri} of a quadruple of coherent vortex formations, but now without introducing artificial hyperviscosity into the equations (see section~\ref{sec:DQMsimulation}).
On the other hand, for non-zero momentum vorticity, our results show that the generic behaviour suggested in \citep{dri} is incorrect: the situation is more complicated yielding either 2, 3, or 4 coherent vortex formations.


In \autoref{sec:pvd}, comparing the motion of the obtained vortex blob formations with point vortex dynamics, we presented a theoretical explanation describing the mechanism for the asymptotic behaviour of the solutions to the Euler equations: the inverse energy cascade continues until two, three or four vortex blobs have been formed, with the number of vortices correlated to the ratio between the magnitude of the momentum and the square root of the enstrophy.
After that, the vortex blob formation is blocked from further mixing by the quasi-periodic motion imposed by the integrability of the point vortex dynamics.
This way, we establish integrability theory of point vortex dynamics together with KAM perturbation theory as the fundamental theory underlying the formation of unsteady but quasi-periodic coherent vortex formations.



As an outlook, the connections to integrability theory could be studied in much more detail, for example, the relation between the regularity of the generic initial data and the qualitative properties of the final state of the system, e.g.\ the size of the vortex blobs.
Perhaps more pressing is to get better statistics for the correlation between $\gamma$ and the long-time behaviour: instead of just 16 simulations, we aim to run 512 or more simulations to collect statistics from.
One could also try with initial vorticity from the more regular Gaussian random field on the Sobolev space $H^s(\Ss^2)$.
Another aspect to investigate is the long-time behaviour of the vorticity on a rotating sphere. 
In the appendix, we present some numerical results indicating that quasi-periodic behaviour can also be reached, but is now more complicated than what can be achieved by point vortex dynamics.
One could also look deeper at the mechanism behind the inverse energy cascade in the quantized equations.
For example, the standard Poisson bracket between two spherical harmonics 
feeds into harmonics with larger wave numbers $l$.
In the quantized bracket, however, high wave number harmonics are fed to lower wave numbers.
This might explain why the inverse energy Cascade works well despite spatial truncation. 
Another benefit of the quantized fluid model is that it is possible to introduce viscosity while still preserving all of the Casimirs.
Indeed, one can add a gradient term of (some approximation of) the entropy functional in such a way that isospectrality is preserved; an example is the \emph{Brockett flow}~\citep{Br1991b} which is known to correspond to a gradient flow of entropy on the space of multivariate Gaussian probability distributions~\citep{Mo2017}.
If viscosity is added to the quantized vorticity equations in such a way, the resulting model can be seen as a mix of computational ideal fluid dynamics (corresponding to the conservative part of the dynamics) and the MRS statistical mechanics model (corresponding to the pure spectral preserving entropy maximizing gradient flow). 



\appendix
\section{Rotating sphere: Rossby-Haurwitz (RH) waves}\label{sec:rossby}
Although the main focus in this paper is on the non-rotating sphere, we like to stress that the numerical method also works well for rotating spheres, of high relevance in geophysical flows.
Indeed, we demonstrate in this appendix how our spatial-time discretization also captures typical features of the quasi-geostrophic equations on a rotating sphere. 
A well known class of exact solutions to the vorticity equation \eqref{euleq2} on a rotating sphere are the Rossby-Haurwitz (RH) waves. 
In terms of spherical harmonics the general formula is
\begin{equation}\label{RH_waves}
	\omega(\phi,\theta,t) = C f + \sum_{m=-l}^l \omega^{lm}Y_{lm}(\phi+2\Omega\alpha_l t,\theta)
\end{equation}
where $\alpha_l=\frac{1}{2}\left(\frac{2C}{l(l+1)}-C+1\right)$, $\omega^{lm}\in\Cc$, $C\in\Rr$ and $l=1,2,\dots$. In particular, for $C=\frac{l(l+1)}{l(l+1)-2}$, we get $\alpha_l=0$ corresponding to \textit{stationary RH waves}. 

That (\ref{RH_waves}) are exact solutions to \eqref{euleq2} depends only on the algebraic properties of the Poisson bracket of the spherical harmonics. 
Indeed, it is not hard to check\footnote{A direct computation, together with the fact that $\exp(-T_{10})\Delta^{-1}_N(A)\exp(T_{10})=\Delta^{-1}_N(\exp(-T_{10})A\exp(T_{10}))$ and $F\varpropto T_{10}$.} that we get an analogous class of exact solutions to (\ref{eq:quant_euler_cor}) in terms of $T^N_{lm}$:

\begin{equation}\label{RH_waves_quant}
W(t) = C\cdot F + \exp(-\alpha_l N^{3/2} F\cdot t)\sum_{m=-l}^l W^{lm}\mathrm{i} T^N_{lm}\exp(\alpha_l N^{3/2} F\cdot t)
\end{equation}
where $\alpha_l=\frac{1}{2}\left(\frac{2C}{l(l+1)}-C+1\right)$, $W^{lm}\in\Cc$, $C\in\Rr$ and $l=1,2,\dots,N$ and $\exp$ is the usual matrix exponential.
We call these solutions \emph{quantized RH waves}.

The stability of RH waves are studied by \citet{Sk2008}.
In essence, they are stable only if they exhibit zonal symmetry.
We have carried out several simulations with our method verifying that the stable exact RH waves correspond to stable quantized RH waves.
We predict that the stability analysis carried out by Skiba can be adopted to the quantized RH waves.

\begin{figure}
	\centering
	\includegraphics[width=0.99\textwidth]{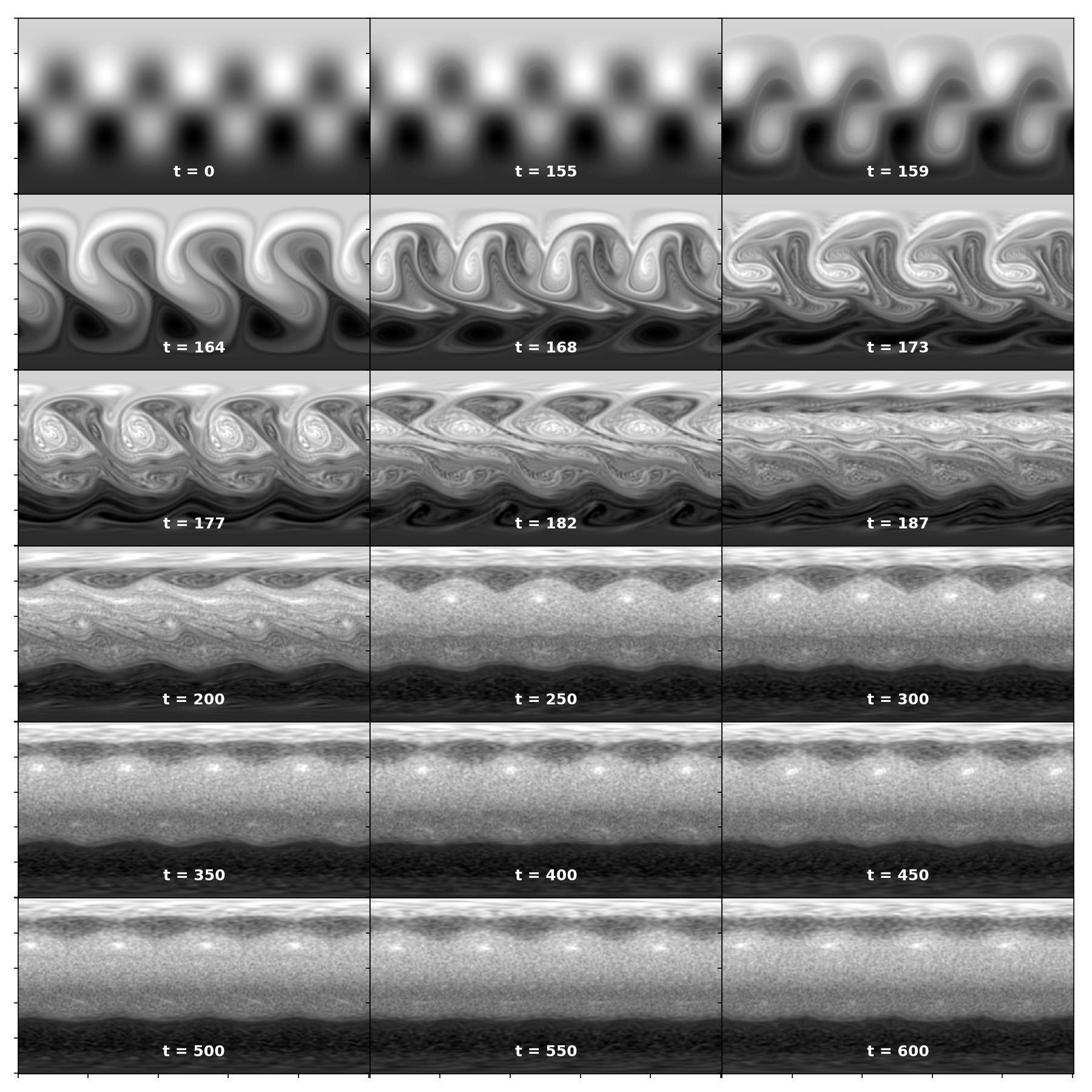}
	\caption{Unsteady quantized RH wave, for the initial conditions as in (\ref{RH_waves_quant}) with parameters \eqref{eq:RH_data}.
	Due to numerical rounding errors, the wave eventually breaks up, goes through an intermediate transition, and then reaches a quasi-periodic asymptotic with sliding zonal vortex belts.
	See also Movie 8 of the supplementary material.
	}\label{fig:unsteady_RH_wave}
\end{figure}

Let us now study the break-up of a non-stable quantized RH wave.
To this end, consider the quantized RH waves with real components
\begin{equation}\label{eq:RH_data}
	C=1, \quad W^{10}=12.9487,\quad W^{54}=W^{5(-4)}=7.7300.
\end{equation}
This wave is non-stable, as it does not have zonal symmetry.
It is also non-stationary.
We use the spatial discretization parameter $N=501$ and the Heun time integration method, with the same non-dimensional parameters as in the previous simulations.
Although the quantized wave is an exact solution to the quantized vorticity equation, due to rounding errors the numerical simulation eventually drift away. 
This can be seen in \autoref{fig:unsteady_RH_wave}.
Up until about $t=155$ the solution remains close to the quantized RH wave.
At $t=159$ it starts to break up in a complicated way.
There is then a transition up until about $t=350$.
After that, the solution settles again on a quasi-periodic asymptotic, but more complicated then in the non-rotating case studied in \autoref{sec:DQMsimulation}.
One can see sliding zonal bands separated by sharp gradients, with `vortex streets' similar in character to those regularly seen on Jupiter \citep{HuMa2007}\footnote{See also \href{https://en.wikipedia.org/wiki/Atmosphere_of_Jupiter}{wikipedia.org/wiki/Atmosphere\_of\_Jupiter}}, see Movie 8 of the supplementary material.
The fluid behaviour shown in \autoref{fig:unsteady_RH_wave} can be found among the regimes described by \citet{NoYo1997}, even though in our simulation smaller vortices inside the alternating jets still persist.


\bibliographystyle{jfm}
\bibliography{biblio1}

\begin{thebibliography}{48}
\expandafter\ifx\csname natexlab\endcsname\relax\def\natexlab#1{#1}\fi
\def\au#1{#1} \def\ed#1{#1} \def\yr#1{#1}\def\at#1{#1}\def\jt#1{\textit{#1}}
  \def\bt#1{#1}\def\bvol#1{\textbf{#1}} \def\vol#1{#1} \def\pg#1{#1}
  \def\publ#1{#1}\def\arxiv#1{#1}\def\org#1{#1}\def\st#1{\textit{#1}}

\bibitem[Abramov \& Majda(2003)]{AbMa2003}
{\sc \au{Abramov, R.~V.} \& \au{Majda, A.~J.}} \yr{2003}  \at{Statistically
  relevant conserved quantities for truncated quasigeostrophic flow}.
  \jt{Proc. Nat. Acad. Sci. (USA)}  \bvol{100}~(7),  \pg{3841--3846}.

\bibitem[Aref(2007{\natexlab{{\em a\/}}})]{are}
{\sc \au{Aref, H.}} \yr{2007{\natexlab{{\em a\/}}}}  \at{Point vortex dynamics:
  A classical mathematics playground}.  \jt{J. Math. Phys.}  \bvol{48, 065401}.

\bibitem[Aref(2007{\natexlab{{\em b\/}}})]{Ha2007}
{\sc \au{Aref, H.}} \yr{2007{\natexlab{{\em b\/}}}}  \at{Point vortex dynamics:
  a classical mathematics playground}.  \jt{J. Math. Phys.}  \bvol{48}~(6),
  \pg{065401, 23}.

\bibitem[Arnold(1966)]{arn}
{\sc \au{Arnold, V.~I.}} \yr{1966}  \at{Sur la g\'eometrie differentielle des
  groupes de lie de dimension infinie et ses applications a l'hydrodynamique
  del fluids parfaits}.  \jt{Ann. Fourier}  \bvol{16}~(1),  \pg{319--361}.

\bibitem[Arnold(1989)]{Ar1989}
{\sc \au{Arnold, V.~I.}} \yr{1989} {\em Mathematical Methods of Classical
  Mechanics\/}, 2nd edn.  \publ{New York: Springer-Verlag}.

\bibitem[Arnold \& Khesin(1998)]{ArKh1998}
{\sc \au{Arnold, V.~I.} \& \au{Khesin, B.~A.}} \yr{1998} {\em Topological
  Methods in Hydrodynamics\/}.  \publ{New York: Springer-Verlag}.

\bibitem[Bogomolov(1977)]{bog}
{\sc \au{Bogomolov, V.A.}} \yr{1977}  \at{{Dynamics of vorticity at sphere}}.
  \jt{Fluid Dyn.}  \bvol{16},  \pg{863--870}.

\bibitem[Bordemann {\em et~al.\/}(1991)Bordemann, Hoppe, Schaller \&
  Schlichenmaier]{bhss}
{\sc \au{Bordemann, M.}, \au{Hoppe, J.}, \au{Schaller, P.} \&
  \au{Schlichenmaier, M.}} \yr{1991}  \at{{$\mathfrak{gl}(\infty)$ and
  geometric quantization}}.  \jt{Comm. Math. Phys.}  \bvol{138}~(2),
  \pg{209--244}.

\bibitem[Bordemann {\em et~al.\/}(1994{\natexlab{{\em a\/}}})Bordemann,
  Meinrenken \& Schlichenmaier]{bms}
{\sc \au{Bordemann, M.}, \au{Meinrenken, E.} \& \au{Schlichenmaier, M.}}
  \yr{1994{\natexlab{{\em a\/}}}}  \at{{Toeplitz Quantization of K{\"a}hler
  Manifolds and $\mathfrak{gl}(N), N\rightarrow\infty$ limits}}.  \jt{Comm.
  Math. Phys.}  \bvol{165}~(2),  \pg{281--296}.

\bibitem[Bordemann {\em et~al.\/}(1994{\natexlab{{\em b\/}}})Bordemann,
  Meinrenken \& Schlichenmaier]{BoMeSc1994}
{\sc \au{Bordemann, M.}, \au{Meinrenken, E.} \& \au{Schlichenmaier, M.}}
  \yr{1994{\natexlab{{\em b\/}}}}  \at{Toeplitz quantization of {K}\"ahler
  manifolds and $\mathfrak{gl}(n), n\to\infty$ limits}.  \jt{Comm. Math. Phys.}
   \bvol{165}~(2),  \pg{281--296}.

\bibitem[Brockett(1991)]{Br1991b}
{\sc \au{Brockett, R.~W.}} \yr{1991}  \at{Dynamical systems that sort lists,
  diagonalize matrices, and solve linear programming problems}.  \jt{Linear
  Algebra Appl.}  \bvol{146},  \pg{79--91}.

\bibitem[Dolzhansky(2012)]{Do2012}
{\sc \au{Dolzhansky, F.~V.}} \yr{2012} {\em Fundamentals of geophysical
  hydrodynamics\/}.  \publ{Springer}.

\bibitem[Dritschel {\em et~al.\/}(2015)Dritschel, Qi \& Marston]{dri}
{\sc \au{Dritschel, D.G.}, \au{Qi, W.} \& \au{Marston, J.B.}} \yr{2015}
  \at{{On the late-time behaviour of a bounded, inviscid two-dimensional
  flow}}.  \jt{J. Fluid Mech.}  \bvol{783},  \pg{1--22}.

\bibitem[Dubinkina \& Frank(2010)]{DuFr2010}
{\sc \au{Dubinkina, S.} \& \au{Frank, J.}} \yr{2010}  \at{Statistical relevance
  of vorticity conservation in the {H}amiltonian particle-mesh method}.  \jt{J.
  Comp. Phys.}  \bvol{229}~(7),  \pg{2634--2648}.

\bibitem[Hairer {\em et~al.\/}(2006)Hairer, Lubich \& Wanner]{HaLuWa2006}
{\sc \au{Hairer, E.}, \au{Lubich, C.} \& \au{Wanner, G.}} \yr{2006} {\em
  Geometric Numerical Integration\/}, 2nd edn.  \publ{Berlin: Springer-Verlag}.

\bibitem[Helmholtz(1858)]{Hel1858}
{\sc \au{Helmholtz, H.}} \yr{1858}  \at{{Über Integrale der hydrodynamischen
  Gleichungen, welche den Wirbelbewegungen entsprechen}}.  \jt{J. Reine Angew.
  Math.}  \bvol{25--55},  \pg{25–55}.

\bibitem[Hoppe(1982)]{hopPhD}
{\sc \au{Hoppe, J.}} \yr{1982}  \at{{Ph.D. thesis MIT Cambridge}} .

\bibitem[Hoppe \& Yau(1998)]{hopyau}
{\sc \au{Hoppe, J.} \& \au{Yau, S.-T.}} \yr{1998}  \at{{Some properties of
  matrix harmonics on S2}}.  \jt{Comm. Math. Phys.}  \bvol{195},  \pg{66--77}.

\bibitem[Humphreys \& Marcus(2007)]{HuMa2007}
{\sc \au{Humphreys, T.} \& \au{Marcus, P.~S.}} \yr{2007}  \at{Vortex street
  dynamics: The selection mechanism for the areas and locations of {J}upiter's
  vortices}.  \jt{J. Atmosph. Sci.}  \bvol{64}~(4),  \pg{1318--1333}.

\bibitem[Keiner \& Potts(2008)]{KePo2008}
{\sc \au{Keiner, J.} \& \au{Potts, D.}} \yr{2008}  \at{Fast evaluation of
  quadrature formulae on the sphere}.  \jt{Math. Comp.}  \bvol{77}~(261),
  \pg{397--419}.

\bibitem[Kraichnan(1967)]{Kr1967}
{\sc \au{Kraichnan, R.~H.}} \yr{1967}  \at{Inertial ranges in two-dimensional
  turbulence}.  \jt{Phys. Fluid.}  \bvol{10}~(7),  \pg{1417--1423}.

\bibitem[Lang \& Schwab(2015)]{LaSc2015}
{\sc \au{Lang, A.} \& \au{Schwab, C.}} \yr{2015}  \at{Isotropic {G}aussian
  random fields on the sphere: Regularity, fast simulation and stochastic
  partial differential equations}.  \jt{Ann. Appl. Probab.}  \bvol{25}~(6),
  \pg{3047--3094}.

\bibitem[Laurent-Polz {\em et~al.\/}(2011)Laurent-Polz, Montaldi \&
  Roberts]{LaMoRo2011}
{\sc \au{Laurent-Polz, F.}, \au{Montaldi, J.} \& \au{Roberts, M.}} \yr{2011}
  \at{Point vortices on the sphere: stability of symmetric relative
  equilibria}.  \jt{J. Geom. Mech.}  \bvol{3}~(4),  \pg{439--486}.

\bibitem[Majda \& Bertozzi(2002)]{MaBe2002}
{\sc \au{Majda, A.~J.} \& \au{Bertozzi, A.~L.}} \yr{2002} {\em Vorticity and
  {I}ncompressible {F}low\/}.  \publ{Cambridge University Press}.

\bibitem[Marsden \& Ratiu(1998)]{marat}
{\sc \au{Marsden, J.E.} \& \au{Ratiu, T.S.}} \yr{1998} {\em {Introduction to
  Mechanics and Symmetry}\/}.  \publ{Springer}.

\bibitem[Marsden \& Weinstein(1983)]{maz}
{\sc \au{Marsden, J.} \& \au{Weinstein, A.}} \yr{1983}  \at{Co-adjoint orbits,
  vortices, and {C}lebsch variables for incompressible fluids}.  \jt{Physica D}
   \bvol{7}~(1--3),  \pg{305--323}.

\bibitem[Marsden {\em et~al.\/}(2007)Marsden, Misio{\l}ek, Ortega, Perlmutter
  \& Ratiu]{MaMiOrRa2007}
{\sc \au{Marsden, J.~E.}, \au{Misio{\l}ek, G.}, \au{Ortega, J.-P.},
  \au{Perlmutter, M.} \& \au{Ratiu, T.~S.}} \yr{2007} {\em Hamiltonian
  Reduction by Stages\/}.  \publ{Berlin: Springer}.

\bibitem[Marsden \& Ratiu(1999)]{MaRa1999}
{\sc \au{Marsden, J.~E.} \& \au{Ratiu, T.~S.}} \yr{1999} {\em Introduction to
  Mechanics and Symmetry\/}, 2nd edn.  \publ{New York: Springer-Verlag}.

\bibitem[McLachlan(1993)]{mcl}
{\sc \au{McLachlan, R.}} \yr{1993}  \at{Explicit {L}ie--{P}oisson integration
  and the {E}uler equations}.  \jt{Phys. Rev. Lett.}  \bvol{71}~(19),
  \pg{3043--3046}.

\bibitem[McLachlan {\em et~al.\/}(2014)McLachlan, Modin \&
  Verdier]{McMoVe2014c}
{\sc \au{McLachlan, R.~I.}, \au{Modin, K.} \& \au{Verdier, O.}} \yr{2014}
  \at{Collective {L}ie--{P}oisson integrators on $\mathbb{R}^3$}.  \jt{IMA J.
  Num. Anal.}  \bvol{35}~(2),  \pg{546--560}.

\bibitem[McLachlan {\em et~al.\/}(2016)McLachlan, Modin \&
  Verdier]{McMoVe2016b}
{\sc \au{McLachlan, R.~I.}, \au{Modin, K.} \& \au{Verdier, O.}} \yr{2016}
  \at{A minimal-variable symplectic integrator on spheres}.  \jt{Math. Comp.}
  \bvol{86}~(307),  \pg{2325--2344}.

\bibitem[Miller(1990)]{Mi1990}
{\sc \au{Miller, J.}} \yr{1990}  \at{Statistical mechanics of {E}uler equations
  in two dimensions}.  \jt{Phys. Rev. Lett.}  \bvol{65},  \pg{2137--2140}.

\bibitem[Modin(2017)]{Mo2017}
{\sc \au{Modin, K.}} \yr{2017}  \at{Geometry of matrix decompositions seen
  through optimal transport and information geometry}.  \jt{J. Geom. Mech.}
  \bvol{9}~(3),  \pg{335--390}.

\bibitem[Modin \& Viviani(2019)]{modviv1}
{\sc \au{Modin, K.} \& \au{Viviani, M.}} \yr{2019}  \at{{L}ie--{P}oisson
  numerical schemes for isospectral flows}.  \jt{Found. Comput. Math.}
  Doi:10.1007/s10208-019-09428-w.

\bibitem[Newton(2016)]{Ne2016}
{\sc \au{Newton, P.~K.}} \yr{2016}  \at{The fate of random initial vorticity
  distributions for two-dimensional {E}uler equations on a sphere}.  \jt{J.
  Fluid Mech.}  \bvol{786},  \pg{1–4}.

\bibitem[Nozawa \& Yoden(1997)]{NoYo1997}
{\sc \au{Nozawa, T.} \& \au{Yoden, S.}} \yr{1997}  \at{Formation of zonal band
  structure in forced two-dimensional turbulence on a rotating sphere}.
  \jt{Phys. Fluids} ~(9),  \pg{2081–93}.

\bibitem[Pedlosky(2013)]{Pe2013}
{\sc \au{Pedlosky, J.}} \yr{2013} {\em Geophysical fluid dynamics\/}.
  \publ{Springer}.

\bibitem[Qi \& Marston(2014)]{QiMa2014}
{\sc \au{Qi, W.} \& \au{Marston, J.~B.}} \yr{2014}  \at{Hyperviscosity and
  statistical equilibria of {E}uler turbulence on the torus and the sphere}.
  \jt{J. Stat. Mech.}  \bvol{2014}~(7),  \pg{P07020}.

\bibitem[Rios \& Straume(2014)]{RiSt2014}
{\sc \au{Rios, P.-M.} \& \au{Straume, E.}} \yr{2014} {\em Symbol
  correspondences for spin systems\/}.  \publ{Springer}.

\bibitem[Robert \& Sommeria(1991)]{RoSo1991}
{\sc \au{Robert, R.} \& \au{Sommeria, J.}} \yr{1991}  \at{Statistical
  equilibrium states for two-dimensional flows}.  \jt{J. Fluid Mech.}
  \bvol{229},  \pg{291–310}.

\bibitem[Sakajo(2007)]{Sa2007}
{\sc \au{Sakajo, T.}} \yr{2007}  \at{Integrable four-vortex motion on sphere
  with zero moment of vorticity}.  \jt{Phys. Fluids} ~(19).

\bibitem[Segre \& Kida(1998)]{SeKi1998}
{\sc \au{Segre, E.} \& \au{Kida, S.}} \yr{1998}  \at{Late states of
  incompressible 2d decaying vorticity fields}.  \jt{Fluid Dyn. Res.}
  \bvol{23}~(2),  \pg{89--112}.

\bibitem[Sevryuk(1995)]{Se1995}
{\sc \au{Sevryuk, M.~B.}} \yr{1995}  \at{K{AM}-stable {H}amiltonians}.  \jt{J.
  Dynam. Control Systems}  \bvol{1}~(3),  \pg{351--366}.

\bibitem[Skiba(2008)]{Sk2008}
{\sc \au{Skiba, Y.~N.}} \yr{2008}  \at{Nonlinear and linear instability of the
  {R}ossby-{H}aurwitz wave}.  \jt{J. Math. Sci.}  \bvol{149}~(6),
  \pg{1708--1725}.

\bibitem[Viviani(May 2019)]{Viv2019}
{\sc \au{Viviani, M.}} \yr{May 2019}  \at{A minimal-variable symplectic method
  for isospectral flows}.  \jt{Arxiv e-prints} .

\bibitem[Zeitlin(1991)]{ze1}
{\sc \au{Zeitlin, V.}} \yr{1991}  \at{{Finite-mode analogues of 2D ideal
  hydrodynamics: Coadjoint orbits and local canonical structure}}.  \jt{Physica
  D}  \bvol{49}~(3),  \pg{353--362}.

\bibitem[Zeitlin(2004)]{ze2}
{\sc \au{Zeitlin, V.}} \yr{2004}  \at{Self-consistent-mode approximation for
  the hydrodynamics of an incompressible fluid on non rotating and rotating
  spheres}.  \jt{Phys. Rev. Lett.}  \bvol{93}~(26),  \pg{353--362}.

\bibitem[Zeitlin(2018)]{Ze2018}
{\sc \au{Zeitlin, V.}} \yr{2018} {\em Geophysical fluid dynamics: understanding
  (almost) everything with rotating shallow water models\/}.  \publ{Oxford
  University Press}.

\end{thebibliography}

\end{document}